\documentclass[11pt]{amsart}
\usepackage[a4paper,margin=29mm]{geometry}
\usepackage{amsmath,amssymb,amsthm,mathtools}
\usepackage[T1]{fontenc}
\usepackage{lmodern}
\usepackage{xcolor}
\usepackage[colorlinks=true,linkcolor=blue,citecolor=blue,urlcolor=blue]{hyperref}
\numberwithin{equation}{section}
\newtheorem{theorem}{Theorem}[section]
\newtheorem{proposition}[theorem]{Proposition}
\newtheorem{lemma}[theorem]{Lemma}
\newtheorem{corollary}[theorem]{Corollary}
\theoremstyle{definition}

\newtheorem{assumption}[theorem]{Assumption}
\theoremstyle{remark}
\newtheorem{remark}[theorem]{Remark}
\newcommand{\R}{\mathbb R}
\newcommand{\T}{\mathbb T}
\newcommand{\Sph}{\mathbb S}

\newcommand{\av}[1]{\langle #1\rangle}
\newcommand{\Tr}{\mathcal T}
\newcommand{\HH}{\mathcal H}
\newcommand{\DD}{\mathcal D}
\newcommand{\mic}{\mathrm{mic}}
\newcommand{\mac}{\mathrm{mac}}
\newcommand{\cP}{\mathcal P}
\newcommand{\norm}[1]{\left\lVert #1\right\rVert}

\DeclareMathOperator{\spanop}{span}
\allowdisplaybreaks[2]
\title[Entropy decay for Boltzmann and BGK]{Constructive Trend to Equilibrium for Low-Regularity Boltzmann and BGK Equations}
\author{Donghyun Lee}
\address{Department of Mathematics, POSTECH}
\email{donlgee@postech.ac.kr}
\date{Draft: October 9, 2026}

\begin{document}
\begin{abstract}
We study the long-time behavior of global mild solutions to the cutoff
Boltzmann and BGK equations on $\T^3$ under a Gaussian lower bound and
a polynomially weighted $L^\infty$ upper bound. The main step is a
finite-time estimate which uses free transport to control the entropy
of the local Maxwellian by microscopic entropy and dissipation, without
assuming Sobolev regularity. For Boltzmann, Villani's entropy-production
inequality and a comparison of collision kernels give polynomial decay
when 
the velocity weight has order $m>7$. For a sufficiently large weight, the large-amplitude theory
then gives exponential decay in the weighted norm and in relative
entropy. For BGK, the microscopic inequality is linear, so entropy
decays exponentially from the finite-time estimate.
{\color{black}
Finally, we establish a density criterion for large-data BGK dynamics
in the polynomially weighted $L^\infty$ class of Perthame and Pulvirenti.
For weight order $m>10$, a uniform upper bound on the local density
guarantees uniform-in-time weighted bounds and exponential convergence
to equilibrium in both relative entropy and the weighted $L^\infty$ norm.
}
\end{abstract}
\maketitle
\tableofcontents

\section{Introduction}\label{sec:intro}

We study the long-time behavior of low-regularity
solutions to the Boltzmann and BGK equations on $\T^3\times\R^3$:
\begin{equation}\label{eq:BEintro}
 \partial_t F+v\cdot\nabla_xF=Q(F,F)
\end{equation}
and
\begin{equation}\label{eq:BGKintro}
 \partial_t F+v\cdot\nabla_xF=M_F-F.
\end{equation}

For the Boltzmann equation, we use the collision parametrization
\[
v'=\frac{v+v_*}{2}+\frac{|v-v_*|}{2}\sigma,\qquad
v_*'=\frac{v+v_*}{2}-\frac{|v-v_*|}{2}\sigma,\qquad\sigma\in\Sph^2,
\]
and
\begin{equation}\label{eq:collision}
	Q(F,F)(v)=\int_{\R^3\times\Sph^2}
	B(v-v_*,\sigma)(F'F_*'-FF_*)d\sigma d v_*.
\end{equation}
We assume 
\begin{equation}\label{eq:kernel}
	B(w,\sigma)=|w|^\gamma b(\widehat w\cdot\sigma),
	\quad 0\le\gamma\le1,\quad b\ge0,\quad
	\beta_b:=\int_{\Sph^2}b(e\cdot\sigma)d\sigma<\infty.
\end{equation}
The value at $w=0$ does not affect the integrals. We use the usual
pre/post-collisional symmetries. 
The finite-time entropy estimate in \eqref{eq:BE-window} requires only the nonnegativity and angular integrability of $b$. For the algebraic entropy decay obtained
using Proposition~\ref{aux:proposition}, we additionally assume that
$b\ge c_b>0$ almost everywhere and that $m>7$. We define the entropy by
\begin{equation} \notag
H(F) := \iint_{\mathbb{T}^3\times\R^3} F\log F dvdx.
\end{equation}
When $F$ is sufficiently smooth and decays sufficiently rapidly at large
velocities, the entropy is nonincreasing in time and satisfies
\begin{equation} \label{H-theorem BE}
\frac{d}{dt}H(F)(t) + \mathcal{D}_B(F)(t) = 0,
\end{equation}
where the entropy dissipation $\mathcal{D}_B$ and the local entropy
dissipation $D_v$ are defined by
\begin{equation}
	\begin{split}
	D_v(f)&=\frac14\int_{\R^3\times\R^3\times\Sph^2}
	B(v-v_*,\sigma)\Psi(f'f_*',ff_*)  d\sigma  d v_*  d v, \\
	\DD_{\rm B}(t)&=\int_{\T^3}D_v(F)  d x.\label{eq:DB}
\end{split}
\end{equation}
where
\begin{equation} \label{eq:Psi}		
	\Psi(u,w)=(u-w)\log(u/w)\ge 0.
\end{equation}
For the BGK equation, $M_F$ is the local Maxwellian determined by the mass,
momentum, and energy of $F(t,x,\cdot)$. We define the macroscopic fields
$(\rho, U, T)$ by
\begin{equation}\label{eq:fields}
	\rho=\int_{\R^3}F  d v,\qquad
	\rho U=\int_{\R^3}vF  d v,\qquad
	3\rho T=\int_{\R^3}|v-U|^2F  d v,
\end{equation}
and the local Maxwellian $M_F$ is given by
\begin{equation}\label{eq:localM}
	M_F(v)=\frac{\rho}{(2\pi T)^{3/2}}
	\exp\!\left(-\frac{|v-U|^2}{2T}\right).
\end{equation}
The entropy for the BGK equation is also nonincreasing in time:
\begin{equation} \label{H-theorem BGK}
	\frac{d}{dt}H(F)(t) + \mathcal{D}_G(F)(t) = 0,
\end{equation}
where
\begin{equation}\label{eq:DG}
	D_{{\rm G},v}(F)=\int_{\R^3}(F-M_F)\log(F/M_F)  d v,
	\qquad   \DD_{\rm G}(t)=\int_{\T^3}D_{{\rm G},v}(F)  d x.
\end{equation} 

For $f,g \geq 0$, we define the relative entropy by
\[
	H(f\mid g) := \iint_{\mathbb{T}^3\times\R^3} \Phi(f\mid g) dvdx = \int_{\mathbb{T}^3} H_v(f\mid g) dx
\]
where
\[
\Phi(f\mid g) :=f\log(f/g)-f+g,
\qquad H_v(f\mid g) :=\int_{\R^3}\Phi(f\mid g)  d v.
\]

We work on $\T^3=(\R/2\pi\mathbb Z)^3$, equipped throughout with the
normalized measure $dx=(2\pi)^{-3} dx_{\rm Leb}$, where $dx_{\rm Leb}$ is
Lebesgue measure on $\mathbb{R}^3$. Thus
$\int_{\T^3}1  d x=1$ and spatial Fourier frequencies are $k\in\mathbb Z^3$.
The velocity measure and the surface measure on $\Sph^2$ are unnormalized.
We assume that, for every relevant time,
\begin{equation}\label{eq:conservation}
	\int_{\T^3\times\R^3}(1,v,|v|^2)F(t,x,v)  d v  d x=(1,0,3).
\end{equation}
Set $\av v=(1+|v|^2)^{1/2}$.  \\

The $H$-theorem provides the basic mechanism for relaxation to
equilibrium, but entropy monotonicity alone does not give a quantitative
rate. A central question, formulated in Cercignani's conjecture, is
whether entropy production controls the relative entropy with respect
to the associated Maxwellian. Villani \cite{V03} established a linear
entropy-production inequality for collision kernels with superquadratic
growth, together with weakened inequalities for physically relevant
kernels under additional a priori bounds. These estimates lead to
quantitative relaxation beyond the perturbative regime, although the
form of the inequality and the resulting rate depend on the collision
kernel and the available bounds on the distribution. We refer to
Desvillettes, Mouhot, and Villani \cite{DMV11} for a survey of the
conjecture, its limitations, and its role in the study of entropy decay.

For the spatially homogeneous Boltzmann equation with cutoff hard
potentials, Mouhot and Villani \cite{MV04} developed a constructive
regularity theory based on the smoothing properties of the gain
operator. Their decomposition of solutions into a uniformly regular
part and an exponentially decaying remainder makes it possible to
apply entropy estimates without assuming high regularity of the
initial datum. Under their assumptions on the kernel and initial data,
this yields convergence faster than any prescribed algebraic rate.
Mouhot \cite{M06} subsequently obtained quantitative exponential
convergence by combining nonlinear relaxation estimates with spectral
estimates near equilibrium. These results illustrate how entropy
methods and the linearized theory can be combined to obtain sharper
decay rates.

Recent developments for singular collision operators have identified
Fisher information as a second Lyapunov functional. Guillen and
Silvestre \cite{GS25} proved that it is nonincreasing along the
spatially homogeneous Landau flow for interaction laws including the
Coulomb potential, and used this estimate to establish global
regularity for nonnegative initial data bounded above by a Maxwellian.
Imbert, Silvestre, and Villani \cite{ISV26} extended this approach to
the spatially homogeneous Boltzmann equation for a large class of
collision kernels, including noncutoff very soft potentials. Their
monotonicity criterion is tied to an integro-differential inequality
of logarithmic Sobolev type on the sphere and yields global smooth
solutions under the stated kernel and initial-data assumptions.
These results complement entropy-production methods by showing how
a monotone information functional can provide the additional
regularity estimates needed for singular interactions.

In the spatially inhomogeneous problem, an additional difficulty is
that collisions relax the distribution towards a local Maxwellian,
whereas convergence to global equilibrium also requires the action of
transport. Desvillettes and Villani \cite{DV05} established convergence
faster than any prescribed algebraic rate under suitable uniform
assumptions on regularity, velocity decay, and strict positivity.
Their work provides a quantitative entropy approach to the full
nonlinear equation, conditional on these a priori bounds.
Gualdani, Mischler, and Mouhot \cite{GMM17} later obtained a
constructive exponential $H$-theorem for hard spheres on the torus,
conditional on suitable uniform smoothness and polynomial moment
estimates, through a factorization method for non-symmetric operators.

The underlying rigidity principle is classical: a local Maxwellian
that solves the free-transport equation on the torus must be a global
Maxwellian. This principle is central to the approach of
Desvillettes and Villani \cite{DV05}, who quantify the departure from
local equilibrium through differential inequalities and geometric
estimates under uniform regularity assumptions. The measurable
version of this rigidity principle and its use in a compactness proof
of qualitative entropy convergence are also discussed in
\cite[Sections 4.1 and 6.2]{HJJLP26}.

Finite-time observability provides another way to quantify the
interaction between transport and dissipation. Tran \cite{Tran13}
established such estimates for conservative transport and used them
to obtain decay for several kinetic models, including the linearized
Boltzmann equation. In particular, his analysis on bounded velocity
sets uses spatial Fourier expansions and estimates of oscillatory
terms. Related space-time coercivity inequalities were developed by
Brigati \cite{Brigati23} for kinetic Fokker--Planck equations: the
space-time norm is controlled by the deviation from the velocity
average and a weak norm of the transport derivative, using averaging
and Lions-type inequalities.

The present transport estimate belongs to this observability approach.
We give an explicit Fourier proof on a bounded velocity ball, with
projection onto the five-dimensional space of collision invariants,
for $L^2$ functions whose transport derivative belongs to $L^2$. We then
apply the estimate to $\log F$: the pointwise envelopes allow us to
control its microscopic component by microscopic entropy and its
transport derivative by entropy dissipation. Together with the
conservation laws, this yields a uniform finite-time estimate for
the macroscopic entropy, without assuming spatial or velocity
Sobolev regularity or small initial relative entropy. It thus gives
a quantitative stability form of the transported-Maxwellian rigidity
principle suited to the nonlinear Boltzmann and BGK equations.

Entropy convergence for low-regularity solutions is particularly
relevant to kinetic boundary problems, where grazing trajectories
and boundary geometry can obstruct smoothness. Discontinuities
can form and propagate in nonconvex domains \cite{Kim11}, while
even in strictly convex domains second-order regularity up to the
boundary can fail \cite{GKTT17}. These phenomena make strong
smoothness assumptions restrictive in the study of relaxation to
equilibrium and motivate entropy estimates formulated directly
for low-regularity solutions. The present analysis develops such
estimates in the periodic setting.

A complementary line of research concerns solutions with large
amplitude but small entropy or other small integral quantities.
Duan, Huang, Wang, and Yang \cite{DHWY17} proved global well-posedness
for initial perturbations that may be large in a weighted
$L^\infty$ norm, subject to smallness of the defect mass, energy,
and entropy together with the $L^1_xL^\infty_v$ norm of the perturbation.
For hard potentials on the torus, their solutions converge
exponentially to equilibrium. Related large-amplitude results in
bounded domains with diffuse reflection were obtained by Duan and
Wang \cite{DW19} under a smallness assumption in $L^2$.
For specular reflection in bounded uniformly convex $C^3$ domains,
Duan, Ko, and Lee \cite{DKL23} established global solutions and
exponential weighted relaxation for cutoff hard potentials,
allowing large initial amplitude under a small relative-entropy
condition.
The polynomially weighted large-amplitude theorem in
\cite{HJJLP26} gives exponential relaxation for an already global
solution with a uniform weighted bound once its relative entropy is
sufficiently small. Thus large amplitude and small relative entropy
play distinct roles in these theories.

For the BGK model, introduced in \cite{BGK54}, the relaxation operator
has a simpler microscopic entropy structure. Perthame and Pulvirenti
\cite{PP93} established weighted $L^\infty$ bounds and uniqueness
for a class of mild solutions, with bounds on each finite time
interval. More recently, Bae, Ko, Lee, and Yun \cite{BKLY26} proved
global existence and exponential weighted relaxation for
large-amplitude initial data with sufficiently small relative
entropy, a suitable polynomial velocity weight, and a positive lower
bound on the free-transported initial density. This provides a BGK
counterpart to the use of entropy smallness in large-amplitude
Boltzmann theory.

Our argument connects these entropy and large-amplitude approaches.
The microscopic Boltzmann estimate starts from Villani's
inequality \cite{V03}. A related comparison with a superquadratic
kernel in \cite[Section 6.1]{HJJLP26} uses H\"older's inequality;
here we separate small and large relative velocities to make the
dependence of the entropy exponent on the polynomial tail explicit.
The finite-time transport estimate applies to both Boltzmann and BGK,
while their different microscopic inequalities yield algebraic and
exponential entropy decay, respectively. No initial entropy
smallness is imposed. For the stronger conclusions in weighted
$L^\infty$, we apply \cite{HJJLP26,BKLY26} only after the entropy
has crossed the required smallness threshold. This last step is
essential: an entropy bound together with a uniform weighted
$L^\infty$ bound does not by itself imply convergence in that
weighted norm. \\

We start with a global mild solution and
assume that, after some time $t_*\ge0$, it is bounded from below by a
Gaussian and from above in a polynomially weighted $L^\infty$ space.
Such bounds are related to quantitative lower-bound results for the
cutoff Boltzmann equation under hydrodynamic assumptions \cite{M05}
and to the conditional large-data theory in \cite{HJJLP26}. Their
propagation is not part of the entropy argument considered here.

Let $\mu(v)=(2\pi)^{-3/2}e^{-|v|^2/2}$ be the global Maxwellian with
the same conserved moments as $F$. The relative entropy splits as
\begin{equation}\label{eq:intro-split}
 \HH(t):=H(F(t)\mid\mu)
 =\underbrace{H(F(t)\mid M_F(t))}_{\HH_{\mic}(t)}
 +\underbrace{H(M_F(t)\mid\mu)}_{\HH_{\mac}(t)}.
\end{equation}
For brevity, we write $h(t,x):=H_v(F(t,x,\cdot)\mid M_F(t,x,\cdot))$.
The collision term controls the first part. It does not directly
control the second: every local Maxwellian has zero collision
dissipation, although its density, velocity, and temperature may vary
with $x$. For instance, if $\rho>0$ and $\int_{\T^3}\rho  d x=1$,
then $F(x,v)=\rho(x)\mu(v)$ has
$\HH_{\mic}=0$ while
$\HH=\int_{\T^3}\rho\log\rho  d x$ can be positive. Transport must
therefore be used before the entropy of $M_F$ can be controlled.

We prove that there exist $T_0,C>0$, independent of the initial time,
such that
\begin{equation}\label{eq:intro-window}
 \int_t^{t+T_0}\HH(s)  d s
 \le C\int_t^{t+T_0}\bigl(\HH_{\mic}(s)+  \DD(s)\bigr)  d s.
\end{equation}
For Boltzmann, we combine this estimate with the microscopic
inequality
\begin{equation}\label{eq:intro-micro}
 D_v(f)\ge c_{\mic}H_v(f\mid M_f)^p,\qquad p>1.
\end{equation}
We establish \eqref{eq:intro-micro} under the additional assumption
$b\ge c_b>0$ almost everywhere and a polynomial upper bound of order
$m>7$. In that case,
\eqref{eq:intro-window} and the entropy inequality yield
\[
 H(F(t)\mid\mu)\le C(1+t-t_*)^{-1/(p-1)}.
\]
Once the entropy is sufficiently small, the large-amplitude
result of \cite{HJJLP26} for polynomial weights gives exponential convergence in a
sufficiently strong weighted $L^\infty$ space. The same rate then
holds for the relative entropy. For BGK,
$D_{{\rm G},v}(F)\ge H_v(F\mid M_F)$, and
\eqref{eq:intro-window} directly gives
\[
 H(F(t)\mid\mu)\le C H(F(t_*)\mid\mu)e^{-\lambda(t-t_*)}.
\]
The large-amplitude result for BGK in \cite{BKLY26} further gives
weighted $L^\infty$ decay when the polynomial weight has order
$m>10$.

The finite-time estimate is based on a local calculation in velocity.
Fix $B=B_1(0)$ and write
\[
 g=\log F,\qquad P=\log M_F=a+b\cdot v+c|v|^2.
\]
On $B$, the Gaussian lower bound turns microscopic entropy into an
$L^2$ estimate of $g-P$. For Boltzmann, the entropy dissipation
controls $\Tr g=Q(F,F)/F$ in $L^2(B)$; for BGK, the same estimate
follows directly from $(M_F-F)/F$. Thus $g$ is close to
\[
 \cP=\spanop\{1,v_1,v_2,v_3,|v|^2\}
\]
and nearly solves the free-transport equation. Here, $\Tr := \partial_t + v\cdot\nabla_x$.
Using Fourier analysis, we show that a polynomial in $\cP$ that nearly follows free
transport must be close to a constant polynomial. Conservation of
mass, momentum, and energy identifies that polynomial with
$\log\mu$. This argument uses the transport graph norm on a bounded
velocity set and does not differentiate the solution in $x$ or $v$.

Section~\ref{sec:setup} states the assumptions and main results.
Section~\ref{sec:transport} proves the transport estimate, and
Section~\ref{sec:macro} applies it to the local Maxwellian.
Sections~\ref{sec:boltzmann} and~\ref{sec:bgk} establish the entropy
decay results. In Section~\ref{sec:application}, we apply the
BGK theorem to the Perthame--Pulvirenti solution class.
Appendix~\ref{sec:aux} proves the microscopic Boltzmann estimate
used above. \\

\section{Main results}\label{sec:setup}

\begin{assumption}
	\label{ass:envelope}
For fixed $t_*\ge0$ and $a_*,b_*,M_0>0$, the solution satisfies
\begin{equation}\label{eq:envelope}
 a_*e^{-b_*|v|^2}\le F(t,x,v)\le M_0\av v^{-m}
 \quad\text{for a.e. }(t,x,v)\in[t_*,\infty)\times\T^3\times\R^3.
\end{equation}
We require $m>5+\gamma$ for the Boltzmann equation and $m>5$ for BGK.
\end{assumption}


For a source $\mathcal C_F$, by a mild solution of $\Tr F=\mathcal C_F$
we mean a measurable function whose characteristic representative satisfies
\begin{equation}\label{eq:mild}
 F(t,y+(t-s)v,v)=F(s,y,v)
 +\int_s^t\mathcal C_F(\tau,y+(\tau-s)v,v)  d\tau
\end{equation}
for almost every $(y,v)$ and all times in each finite interval, with
the locally integrable representative understood. The source is
$Q(F,F)$ or $M_F-F$. The pointwise and moment bounds used below ensure
that this integral is well defined. We use the absolutely continuous representative when
writing the entropy identity. \\

\begin{lemma} 
	\label{lem:macro-bounds}
Assumption~\ref{ass:envelope} with $m>5$ implies
\begin{equation}\label{eq:macro-bounds}
 0<\rho_-\le\rho\le\rho_+,\qquad
 |U|\le U_+,\qquad 0<T_-\le T\le T_+.
\end{equation}
The constants depend only on the two envelopes. In particular,
$F$ and $M_F$ are bounded above and below by positive constants on
every fixed velocity ball $B_R$.
\end{lemma}
\begin{proof}
The bounds on $F$ give
\begin{align*}
 \rho_-&=a_*\int e^{-b_*|v|^2}  d v,&
 \rho_+&=M_0\int\av v^{-m}  d v,\\
 E_+&=M_0\int|v|^2\av v^{-m}  d v,&
 U_+&=(E_+/\rho_-)^{1/2},\qquad T_+=E_+/(3\rho_-).
\end{align*}
The lower bound on the temperature follows from the symmetry of the Gaussian:
\[
 3\rho T\ge a_*\int|v-U|^2e^{-b_*|v|^2}  d v
 \ge a_*\int|v|^2e^{-b_*|v|^2}  d v.
\]
Thus $T_-=(a_*/(3\rho_+))\int|v|^2e^{-b_*|v|^2}  d v>0$.
On $B_R$, the corresponding bounds for the Maxwellian are
\[
 \frac{\rho_-}{(2\pi T_+)^{3/2}}
 e^{-(R+U_+)^2/(2T_-)}\le M_F
 \le\frac{\rho_+}{(2\pi T_-)^{3/2}}.
\]
Together with \eqref{eq:envelope}, these yield the asserted bounds
on $B_R$. They also give $M_F(v)\le C e^{-c|v|^2}$ for all $v$.
\end{proof}

The envelope bounds also give
\begin{equation}\label{eq:log-bound}
 |\log F|+|\log M_F|+|\log\mu|\le C\av v^2.
\end{equation}
Hence $h(t,x):=H_v(F(t,x,\cdot)\mid M_F(t,x,\cdot))$ and $\HH(t)$ are bounded uniformly in time.
Since $\log(M_F/\mu)$ is a linear combination of the collision
invariants, local moment matching gives
\[
 \int(F-M_F)\log(M_F/\mu)  d v=0.
\]
Expanding the three entropy densities proves the exact decomposition
\eqref{eq:intro-split}.

\begin{theorem} 
\label{thm:boltzmann}
Let $F$ be a global mild solution of \eqref{eq:BEintro} satisfying
\eqref{eq:conservation}, \eqref{eq:kernel}, and
Assumption~\ref{ass:envelope} with $m>5+\gamma$.
Then there exist $T_0,C_{\rm B}>0$, depending on the torus, the
envelopes, $\gamma$, and $\beta_b$, such that
\begin{equation}\label{eq:BE-window}
 \int_t^{t+T_0}\HH(s)  d s
 \le C_{\rm B}\int_t^{t+T_0}
                  (\HH_{\mic}(s)+  \DD_{\rm B}(s))  d s,
 \qquad t\ge t_*.
\end{equation}
The finite-time estimate \eqref{eq:BE-window} does not require a
positive lower bound on the angular factor $b$.
For the algebraic decay conclusion, assume in addition that there
exists $c_b>0$ such that $b\ge c_b$ almost everywhere, and that $m>7$.
Set
\[
 \varepsilon_*:=\max\left\{\frac{\gamma}{3},
                    \frac{2-\gamma}{m-7}\right\},
 \qquad p:=1+\varepsilon_*,\qquad r:=\varepsilon_*^{-1}.
\]
By Proposition~\ref{aux:proposition}, there is $c_{\mic}>0$,
independent of $(t,x)$, such that
\begin{equation}\label{eq:MI}
 D_v(F(t,x,\cdot))\ge c_{\mic}h(t,x)^p.
\end{equation}
Consequently,
\begin{equation}\label{eq:BE-decay}
 \HH(t)\le C(1+t-t_*)^{-r},\qquad
 \norm{F(t)-\mu}_{L^1_{x,v}}\le\sqrt{2C}(1+t-t_*)^{-r/2}.
\end{equation}
Here $C$ depends on the torus, the envelope and collision-kernel
parameters, and $\HH(t_*)$.
\end{theorem}

\begin{corollary}
\label{coro:boltzmann}
Assume the hypotheses of Theorem~\ref{thm:boltzmann}, with
$0<\gamma\le1$, $0<{\color{black}c_b}\le b\le b_1<\infty$ almost everywhere, and
\begin{equation}\label{eq:BE-weight-range}
 m>\max\{7,\beta_\infty\},
\end{equation}
where $\beta_\infty$ is the polynomial-weight threshold in
\cite[Theorem 7.26]{HJJLP26} for this collision kernel. More
explicitly, in the notation of that theorem,
\begin{equation*}
	\begin{split}
		\beta_\infty &=\max\left\{1+\gamma+\frac{4b_1}{{\color{black}c_b}},
		\left(\frac{4C_1}{\nu_0}\right)^{4/(\gamma+3)}\right\},  \\
		\nu_0 &=\inf_{v\in\R^3}\frac{1}{\av v^\gamma}
		\int_{\R^3\times\Sph^2}
		|v-v_*|^\gamma b(\widehat{v-v_*}\cdot\sigma)
		\mu(v_*)  d\sigma  d v_*>0,		
	\end{split}
\end{equation*}
with $C_1$ the explicit kernel constant specified there.
Then there are $C_\infty,C_H,\lambda_\infty>0$ such that
\begin{equation}\label{eq:BE-exp-decay}
 \norm{\av v^m(F(t)-\mu)}_{L^\infty_{x,v}}
 \le C_\infty e^{-\lambda_\infty(t-t_*)},\qquad
 \HH(t)\le C_H e^{-\lambda_\infty(t-t_*)},\quad t\ge t_*.
\end{equation}
The constants may depend on the fixed envelope and kernel parameters
and on $\HH(t_*)$; no smallness is imposed at $t_*$.
The weight in \eqref{eq:BE-exp-decay} is the same as in the upper envelope.
\end{corollary}

\begin{theorem}
	\label{thm:bgk}
Let $F$ be a global mild solution of \eqref{eq:BGKintro} satisfying
\eqref{eq:conservation} and Assumption~\ref{ass:envelope} with $m>5$.
For the same geometric time scale $T_0$, there is a computable $K_{\rm G}>0$
depending only on the envelopes and the torus such that
\begin{equation}\label{eq:BGK-window-main}
 \int_t^{t+T_0}\HH(s)  d s
 \le K_{\rm G}\int_t^{t+T_0}  \DD_{\rm G}(s)  d s.
\end{equation}
Set $q=K_{\rm G}/(K_{\rm G}+T_0)\in(0,1)$ and
$\lambda=T_0^{-1}\log(1/q)$. Then
\begin{equation}\label{eq:BGK-decay-main}
 \HH(t)\le q^{-1}\HH(t_*)e^{-\lambda(t-t_*)},\qquad t\ge t_*.
\end{equation}
In particular,
$\norm{F(t)-\mu}_{L^1_{x,v}}\le
\sqrt{2q^{-1}\HH(t_*)}\,e^{-\lambda(t-t_*)/2}$.
If $m>10$, the large-amplitude result \cite[Theorem 1.1]{BKLY26}
also gives constants $C_\infty,\kappa>0$ such that
\begin{equation}\label{eq:BGK-decay-norm}
 \norm{\av v^m(F(t)-\mu)}_{L^\infty_{x,v}}
 \le C_\infty e^{-\kappa(t-t_*)},\qquad t\ge t_*.
\end{equation}
By Lemma~\ref{lem:weighted-to-entropy}, the same norm estimate
gives $\HH(t)\le C_H e^{-\kappa(t-t_*)}$. Together with
\eqref{eq:BGK-decay-main}, this allows the larger of the two
rates $\lambda$ and $\kappa$ to be used in the entropy bound.
\end{theorem}

Unlike the Boltzmann equation, the BGK equation has a linear microscopic
inequality that does not require a separate comparison of kernels. Neither theorem assumes
that the entropy is initially small.

Perthame and Pulvirenti \cite{PP93} constructed a unique BGK mild
solution with a polynomially weighted $L^\infty$ bound on every finite
time interval, assuming a positive initial lower velocity profile.
If its local density is also bounded above for all time, we obtain
the two envelopes required by Theorem~\ref{thm:bgk}.

\begin{theorem}
	\label{thm:PP-application}
	Let $\beta>10$ and let $F_0\ge0$ satisfy
	\begin{equation}\label{eq:app-data}
		A_{\beta,0}:=\norm{\av v^\beta F_0}_{L^\infty_{x,v}}<\infty,
		\qquad F_0(x,v)\ge\varphi(v)\ge0,\qquad
		\int_{\R^3}\varphi(v)  d v>0,
	\end{equation}
	together with the normalized mass, momentum, and energy in
	\eqref{eq:conservation}. Let $F$ be the unique global mild solution
	of the unit-frequency BGK equation constructed in \cite{PP93}.
	Assume in addition that
	\begin{equation}\label{eq:app-rho-upper}
		\sup_{t\ge0,\;x\in\T^3}\rho(t,x)
		\le A_\rho<\infty.
	\end{equation}
	Then there exist $t_1<\infty$ and constants $a,b,M_\beta>0$ such that,
	for a.e. $t\ge t_1$ and $(x,v)$,
	\begin{equation}\label{eq:app-envelopes}
		a e^{-b|v|^2}\le F(t,x,v)
		\le M_\beta\av v^{-\beta}.
	\end{equation}
	Consequently, Theorem~\ref{thm:bgk} gives both exponential entropy
	decay and weighted relaxation. In particular, there exist a finite
	time $t_{\rm ent}\ge t_1$ and constants
	$C_\infty,C_H,\kappa>0$ such that
	\begin{equation}\label{eq:app-weighted-decay}
		\norm{\av v^\beta(F(t)-\mu)}_{L^\infty_{x,v}}
		\le C_\infty e^{-\kappa(t-t_{\rm ent})},
		\qquad
		\HH(t)\le C_H e^{-\kappa(t-t_{\rm ent})},
		\quad t\ge t_{\rm ent}.
	\end{equation}
	The estimates in the proof make the dependence on $A_\rho$ explicit.
	The density bound \eqref{eq:app-rho-upper} is not a conclusion
	of the Perthame--Pulvirenti theorem.  \\
\end{theorem}

\section{A finite-time transport estimate}\label{sec:transport}

Let $B=B_1(0)$, and let $\Pi$ denote orthogonal projection in $L^2(B)$
onto $\cP=\spanop\{1,v_1,v_2,v_3,|v|^2\}$. It acts only in velocity.
For complex Fourier coefficients we use the complexification of this
space.

\begin{lemma}[Transport coercivity]\label{lem:transport}
There is $T_0>0$, determined by $B$ and the torus, such that for every interval $I$ of length $T_0$ and every
real-valued $g$ with $g,\Tr g\in L^2(I\times\T^3\times B)$,
\begin{equation}\label{eq:transport}
 \inf_{P_0\in\cP}\norm{g-P_0}_2^2
 \le C\left(\norm{(I-\Pi)g}_2^2+T_0^2\norm{\Tr g}_2^2\right),
\end{equation}
for some generic constant $C>0$. We interpret $\Tr g$ in distributions. The norms are taken over
$I\times\T^3\times B$, and $P_0(v)$ is independent of $t,x$.
\end{lemma}

\begin{proof}
By time translation assume $I=(0,L)$, where $L=T_0$ will be chosen.
Let
\[
 g_k(t,v)=\int_{\T^3}g(t,x,v)e^{-ik\cdot x}  d x,
 \qquad r_k=(\Tr g)_k.
\]
The spatial measure is normalized, so no volume factor appears in
Parseval's identity. For every $k\in\mathbb Z^3$,
\[
 \partial_tg_k+i(k\cdot v)g_k=r_k
\]
in distributions. Since $B$ is bounded, multiplication by $k\cdot v$
is bounded for each fixed $k$. Alternatively, multiplication by
$e^{itk\cdot v}$ reduces the equation to an ordinary weak time
derivative. Thus $e^{itk\cdot v}g_k\in H^1(I;L^2(B))$.
We use the corresponding continuous $L^2(B)$-valued
representative of $g_k$. Duhamel principle gives
\begin{equation}\label{eq:mode-duhamel}
 g_k(t)=U_k(t-s)g_k(s)+\int_s^tU_k(t-\tau)r_k(\tau)  d\tau,
 \qquad U_k(\tau)f=e^{-i\tau k\cdot v}f,
\end{equation}
with the usual oriented integral for $t<s$. Each $U_k$ is unitary.

Choose a real orthonormal polynomial basis $\psi_1,\ldots,\psi_5$ of
$\cP$. For $f\in L^2(B;\mathbb C)$, write
\[
 \Pi f=\sum_{j=1}^5z_j\psi_j,\qquad
 z_j=\int_B f(v)\psi_j(v)\,d v.
\]
Applying $U_k(\tau)$ and then projecting once more onto $\cP$ gives
\[
 \Pi U_k(\tau)\Pi f
 =\sum_{i=1}^5\left(\sum_{j=1}^5A_{ij}(\tau k)z_j\right)\psi_i,
\]
where
\[
 A_{ij}(\xi)=\int_B\psi_i(v)\psi_j(v)e^{-i\xi\cdot v}  d v,
 \qquad\xi=\tau k.
\]
Indeed, $A_{ij}(\tau k)$ is the coefficient of $\psi_i$ in
$\Pi U_k(\tau)\psi_j$. 

For $\xi\ne0$, put $q_{ij}=\psi_i\psi_j$. Since
$e^{-i\xi\cdot v}=i|\xi|^{-2}\xi\cdot\nabla_v e^{-i\xi\cdot v}$,
integration by parts on $B$ gives
\begin{align*}
 A_{ij}(\xi)
 &=-\frac{i}{|\xi|^2}\int_B
       (\xi\cdot\nabla q_{ij})e^{-i\xi\cdot v}  d v\\
 &\quad+\frac{i}{|\xi|^2}\int_{\partial B}
       q_{ij}(\xi\cdot n)e^{-i\xi\cdot v}  d S_v.
\end{align*}
The volume and boundary terms are both controlled by
$|\xi|^{-1}$. In particular,
\[
 |A_{ij}(\xi)|\le\frac{N_{ij}}{|\xi|},\qquad
 N_{ij}=\norm{\nabla(\psi_i\psi_j)}_{L^1(B)}
       +\norm{\psi_i\psi_j}_{L^1(\partial B)}.
\]
The boundary contribution is retained in $N_{ij}$.
Set $C_B=(\sum_{i,j}N_{ij}^2)^{1/2}$.
Orthonormality of the basis and Cauchy--Schwarz in the finite sum give,
for $\tau k\ne0$,
\begin{align*}
 \norm{\Pi U_k(\tau)\Pi f}_{L^2(B)}^2
 &=\sum_{i=1}^5\left|\sum_{j=1}^5A_{ij}(\tau k)z_j\right|^2\\
 &\le\left(\sum_{i,j=1}^5|A_{ij}(\tau k)|^2\right)
       \sum_{j=1}^5|z_j|^2
 \le\frac{C_B^2}{|k|^2|\tau|^2}\norm{f}_{L^2(B)}^2.
\end{align*}
Here we used $\sum_j|z_j|^2=\norm{\Pi f}_{L^2(B)}^2
\le\norm{f}_{L^2(B)}^2$.
Every nonzero frequency $k\in\mathbb Z^3$ satisfies $|k|\ge1$.
Therefore, with the choice
\begin{equation}\label{eq:T0-choice}
 \tau_0=\max\{1,2C_B\},\qquad L=T_0=3\tau_0,
\end{equation}
we have
\[
 \norm{\Pi U_k(\tau)\Pi}_{L^2(B)\to L^2(B)}
 \le\frac{C_B}{|k||\tau|}
 \le\frac{C_B}{\tau_0}\le\frac12
 \qquad(k\ne0,\ |\tau|\ge\tau_0).
\]

In particular,
\begin{equation}\label{eq:compressed-mixing}
 \norm{\Pi U_k(\tau)\Pi}\le\tfrac12
 \quad(k\ne0,\ |\tau|\ge\tau_0).
\end{equation}
The estimate is uniform in $k$. The choice $L=3\tau_0$ leaves two
subintervals, $(0,\tau_0)$ and $(2\tau_0,3\tau_0)$, separated by
at least $\tau_0$. Once the polynomial basis is fixed, this choice
depends only on the ball and the fixed torus normalization. \\

Fix $k\ne0$ and write $p_k=\Pi g_k$, $h_k=(I-\Pi)g_k$. 
Take $s\in(0,\tau_0)$ and $t\in(2\tau_0,3\tau_0)$.
Applying $\Pi$ to \eqref{eq:mode-duhamel}, decomposing
$g_k(s)=p_k(s)+h_k(s)$, and using $\Pi p_k(s)=p_k(s)$, we obtain
\[
 p_k(t)=\Pi U_k(t-s)\Pi p_k(s)
       +\Pi U_k(t-s)h_k(s)
       +\int_s^t\Pi U_k(t-\tau)r_k(\tau)\,d\tau.
\]
Since $t-s\ge\tau_0$, the first term is controlled by
\eqref{eq:compressed-mixing}. For the second term, $\Pi$ is a
contraction and $U_k$ is unitary. The same properties and
Cauchy--Schwarz in time give
\begin{align*}
 \left\lVert\int_s^t\Pi U_k(t-\tau)r_k(\tau)\,d\tau
       \right\rVert_{L^2(B)}
 &\le\int_s^t\norm{r_k(\tau)}_{L^2(B)}\,d\tau\\
 &\le\sqrt{t-s}\,\norm{r_k}_{L^2(I\times B)}
 \le\sqrt L\,\norm{r_k}_{L^2(I\times B)}.
\end{align*}
The same argument with $s$ and $t$ interchanged applies to the
backward Duhamel formula: the source integral changes orientation,
while \eqref{eq:compressed-mixing} depends only on the absolute
time difference. Thus, with the fixed-time norms taken over $B$,
\begin{align*}
 \norm{p_k(t)}_2&\le\tfrac12\norm{p_k(s)}_2
                       +\norm{h_k(s)}_2+\sqrt L \norm{r_k}_{L^2(I\times B)},\\
 \norm{p_k(s)}_2&\le\tfrac12\norm{p_k(t)}_2
                       +\norm{h_k(t)}_2+\sqrt L \norm{r_k}_{L^2(I\times B)}.
\end{align*}

Substitution yields
\[
 \norm{p_k(s)}_2\le\tfrac23\norm{h_k(s)}_2
              +\tfrac43\norm{h_k(t)}_2+2\sqrt L \norm{r_k}_{L^2(I\times B)}.
\]
Using $(a+b+c)^2\le3(a^2+b^2+c^2)$ yields
\[
 \norm{p_k(s)}_2^2
 \le\frac43\norm{h_k(s)}_2^2
    +\frac{16}{3}\norm{h_k(t)}_2^2
    +12L\norm{r_k}_{L^2(I\times B)}^2.
\]
Integrate over $(s,t)\in(0,\tau_0)\times(2\tau_0,3\tau_0)$
and divide by $\tau_0$. The two time intervals are disjoint,
and $12L\tau_0=4L^2$ because $L=3\tau_0$.
This gives the first estimate below. The second follows from
$\norm{g_k(s)}_2^2=\norm{p_k(s)}_2^2+\norm{h_k(s)}_2^2$,
which holds by orthogonality of $\Pi$.  
\begin{equation}
\begin{split}
 \int_0^{\tau_0}\norm{p_k(s)}_2^2  d s
 &\le\tfrac{16}{3}\norm{h_k}_{L^2(I\times B)}^2+4L^2 \norm{r_k}_{L^2(I\times B)}^2, \\
 \int_0^{\tau_0}\norm{g_k(s)}_2^2  d s
&\le\tfrac{19}{3}\norm{h_k}_{L^2(I\times B)}^2+4L^2 \norm{r_k}_{L^2(I\times B)}^2.
\end{split}
\end{equation}
For any $t\in I$ and $s\in(0,\tau_0)$, unitarity in
\eqref{eq:mode-duhamel} gives
\begin{equation}
	\begin{split}
		\norm{g_k(t)}_2^2 &\leq \big( \norm{g_k(s)}_2 + \sqrt{L}\norm{r_k}_{L^2(I\times B)} \big)^2 \\
			&\le2\norm{g_k(s)}_2^2+2L \norm{r_k}_{L^2(I\times B)}^2.
		\end{split}
	\end{equation}
Averaging in $s$ and integrating in $t$ gives
\[
 \norm{g_k}_{L^2(I\times B)}^2
 \le\frac{2L}{\tau_0}\int_0^{\tau_0}\norm{g_k(s)}_2^2\,d s
       +2L^2\norm{r_k}_{L^2(I\times B)}^2.
\]
Since $L/\tau_0=3$, the preceding estimate implies 
\begin{equation}\label{eq:nonzero-mode}
 \norm{g_k}_{L^2(I\times B)}^2
 \le C(\norm{h_k}_{L^2(I\times B)}^2+L^2 \norm{r_k}_{L^2(I\times B)}^2),
\end{equation}
for some generic $C>0$. \\

For $k=0$, the equation reduces to $\partial_tg_0=r_0$, so
$g_0\in H^1(I;L^2(B))$. Set
\[
 \bar g_0=\frac1L\int_Ig_0(t)\,d t,\qquad
 P_0=\Pi\bar g_0,\qquad h_0=(I-\Pi)g_0.
\]
To estimate the fluctuation around the time average, write
\[
 g_0(t)-\bar g_0=\frac1L\int_I\bigl(g_0(t)-g_0(s)\bigr)\,d s,
 \qquad
 g_0(t)-g_0(s)=\int_s^t r_0(\tau)\,d\tau.
\]
We easily get
\begin{align*}
 \norm{g_0(t)-\bar g_0}_{L^2(B)}^2
 &\le\frac1L\int_I\norm{g_0(t)-g_0(s)}_{L^2(B)}^2\,d s\\
 &\le\frac1L\int_I|t-s|\,d s\,
               \norm{r_0}_{L^2(I\times B)}^2
 \le L\norm{r_0}_{L^2(I\times B)}^2.
\end{align*}
The bounded linear
operator $I-\Pi$ commutes with the time integral, and Jensen gives
\[
 L\norm{(I-\Pi)\bar g_0}_{L^2(B)}^2
 =L\left\lVert\frac1L\int_Ih_0(t)\,d t\right\rVert_{L^2(B)}^2
 \le\norm{h_0}_{L^2(I\times B)}^2.
\]
Finally, decompose $g_0-P_0=(g_0-\bar g_0)+(I-\Pi)\bar g_0$.
The second summand is independent of time, and
$\int_I(g_0(t)-\bar g_0)\,d t=0$ in $L^2(B)$, so the cross term
vanishes after integration in time. Combining these observations,
\begin{align*}
 \norm{g_0-P_0}_{L^2(I\times B)}^2
 &=\norm{g_0-\bar g_0}_{L^2(I\times B)}^2
       +L\norm{(I-\Pi)\bar g_0}_{L^2(B)}^2\\
 &\le L^2\norm{r_0}_{L^2(I\times B)}^2
       +\norm{h_0}_{L^2(I\times B)}^2.
\end{align*}
Since $P_0$ is independent of $x$, subtracting it changes only the
zero spatial Fourier mode. Parseval's identity therefore gives
\[
 \norm{g-P_0}_{L^2(I\times\T^3\times B)}^2
 =\norm{g_0-P_0}_{L^2(I\times B)}^2
   +\sum_{k\ne0}\norm{g_k}_{L^2(I\times B)}^2.
\]

Now, we are ready to prove \eqref{eq:transport}. Applying the zero-mode estimate and \eqref{eq:nonzero-mode}, and
enlarging the constant so that $C\ge1$, we obtain
\begin{align*}
 \norm{g-P_0}_{L^2(I\times\T^3\times B)}^2
 &\le \norm{h_0}_{L^2(I\times B)}^2
       +L^2\norm{r_0}_{L^2(I\times B)}^2\\
 &\quad+C\sum_{k\ne0}\left(
       \norm{h_k}_{L^2(I\times B)}^2
       +L^2\norm{r_k}_{L^2(I\times B)}^2\right)\\
 &\le C\sum_{k\in\mathbb Z^3}\left(
       \norm{h_k}_{L^2(I\times B)}^2
       +L^2\norm{r_k}_{L^2(I\times B)}^2\right).
\end{align*}
The same constant works for every $k$, which allows it to be taken
outside the sum. Since $\Pi$ acts only in velocity and is bounded
on $L^2(B)$, it commutes with the spatial Fourier transform. Hence
$h_k=((I-\Pi)g)_k$, while $r_k=(\Tr g)_k$ by definition.
Parseval in $x$, followed by integration in $(t,v)$, gives
\begin{align*}
 \sum_{k\in\mathbb Z^3}\norm{h_k}_{L^2(I\times B)}^2
 &=\norm{(I-\Pi)g}_{L^2(I\times\T^3\times B)}^2,\\
 \sum_{k\in\mathbb Z^3}\norm{r_k}_{L^2(I\times B)}^2
 &=\norm{\Tr g}_{L^2(I\times\T^3\times B)}^2.
\end{align*}
The interchange of the sums and integrals is justified by Tonelli's
theorem, since all terms are nonnegative. These sums are finite by
the hypotheses on $g$ and $\Tr g$. Recalling $L=T_0$, we conclude
\[
 \norm{g-P_0}_2^2
 \le C\left(\norm{(I-\Pi)g}_2^2+T_0^2\norm{\Tr g}_2^2\right),
\]
where the norms are now over $I\times\T^3\times B$.
To check that this $P_0$ is an admissible real polynomial, note that
$g_0(t,v)=\int_{\T^3}g(t,x,v)\,d x$ is real-valued, and so is
its time average $\bar g_0$. In the chosen real orthonormal basis,
\[
 P_0(v)=\sum_{j=1}^5\left(\int_B\bar g_0(w)\psi_j(w)\,d w\right)
             \psi_j(v).
\]
Each coefficient is real and independent of $t,x$, so $P_0\in\cP$.
Finally,
\[
 \inf_{P\in\cP}\norm{g-P}_2^2\le\norm{g-P_0}_2^2,
\]
which proves \eqref{eq:transport}. No claim that this particular
$P_0$ attains the infimum is needed. 
\end{proof}

\begin{remark}
The estimate uses only $g$ and $\Tr g$ in $L^2$. In particular,
no derivative in $v$ is needed. The decay of
$\Pi U_k(\tau)\Pi$ at a fixed time separation holds for every
$k\ne0$, which is why the same $T_0$ works for all spatial modes. \\
\end{remark}

\section{Control of the macroscopic entropy}\label{sec:macro}

We now apply the transport estimate to the local Maxwellian. Write
\begin{equation}\label{eq:coefficients}
 \phi(v)=(1,v_1,v_2,v_3,|v|^2)^{\mathsf T},\qquad
 P=\log M_F=z\cdot\phi,
\end{equation}
where $z=(a,b,c)\in\R^5$ and
\begin{equation}\label{eq:z-formula}
 a=\log\rho-\tfrac32\log(2\pi T)-\frac{|U|^2}{2T},
 \qquad b=\frac UT,\qquad c=-\frac1{2T}.
\end{equation}
The equilibrium coefficients are
$z_\mu=(-\tfrac32\log(2\pi),0,-\tfrac12)$.

\begin{lemma} 
\label{lem:log-distance}
Under \eqref{eq:envelope} with $m>5$ there is $C_{\log}<\infty$ such
that, for a.e. $t$,
\begin{equation}\label{eq:log-distance}
 \norm{\log F-\log M_F}_{L^2(\T^3\times B)}^2
 \le C_{\log}\HH_{\mic}(t).
\end{equation}
\end{lemma}
\begin{proof}
Choose $0<\ell\le u<\infty$ bounding both $F$ and $M_F$ on $B$,
as in Lemma~\ref{lem:macro-bounds}. Recall $\Phi(r\mid s) :=r\log(r/s)-r+s$. For positive $r,s$,
\[
 \Phi(r\mid s)\ge\frac{(r-s)^2}{2\max\{r,s\}}.
\]
This follows from Taylor's formula for $r\mapsto\Phi(r\mid s)$,
whose second derivative is $1/r$. The mean value theorem for the
logarithm then gives
\[
 |\log r-\log s|^2\le\ell^{-2}|r-s|^2
 \le2u\ell^{-2}\Phi(r\mid s),\qquad r,s\in[\ell,u].
\]
Integrate on the ball and enlarge the nonnegative entropy integral
to all velocities.
\end{proof}

Recall $\phi$ given in \eqref{eq:coefficients}.
Direct computation gives
\[
	\int_B |(z-z_0)\cdot\phi(v)|^2  d v
	= (z-z_0)^{\mathsf T}G_B(z-z_0).
\]
where the matrix $G_B$ is defined by 
\begin{equation}\label{eq:Gram}
 G_B=\int_B\phi\phi^{\mathsf T}  d v
 =\begin{pmatrix}
 4\pi/3&0&0&0&4\pi/5\\
 0&4\pi/15&0&0&0\\
 0&0&4\pi/15&0&0\\
 0&0&0&4\pi/15&0\\
 4\pi/5&0&0&0&4\pi/7
 \end{pmatrix}
\end{equation}
which is positive definite: a polynomial $\xi\cdot\phi$ that vanishes
on the open ball vanishes identically. Denote its smallest eigenvalue
by $\lambda_B>0$. 
Consequently,
\begin{equation}\label{eq:Gram-equivalence}
 \lambda_B|z-z_0|^2
 \le\int_B|(z-z_0)\cdot\phi(v)|^2  d v
 \le\norm{G_B}|z-z_0|^2.
\end{equation}

\begin{lemma} 
\label{lem:anchoring}
Assume \eqref{eq:macro-bounds} and \eqref{eq:conservation}. On any
finite time interval $I$ of length $L$, let
$\bar z=L^{-1}\int_{I\times\T^3}z  d x  d t$.
There is a constant depending only on the macroscopic bounds such that
\begin{equation}\label{eq:anchoring}
 \norm{z-z_\mu}_{L^2(I\times\T^3)}^2
 \le C_{\rm avg}\norm{z-\bar z}_{L^2(I\times\T^3)}^2.
\end{equation}
Moreover,
\begin{equation}\label{eq:macro-equivalence}
\begin{aligned}
 c_{\rm mac}\int_{\T^3}
 (|\rho-1|^2+|U|^2+|T-1|^2)  d x
 &\le \HH_{\mac}(t)\\
 &\le C_{\rm mac}\int_{\T^3}
 (|\rho-1|^2+|U|^2+|T-1|^2)  d x\\
 &\le C'_{\rm mac}\int_{\T^3}|z-z_\mu|^2  d x.
\end{aligned}
\end{equation}
\end{lemma}
\begin{proof}
The bounds \eqref{eq:macro-bounds} place $z$ in a fixed compact set
with $c\le-1/(2T_+)<0$. Choose a compact convex box $K\subset\{c<0\}$
containing this set and $z_\mu$. Then $\bar z\in K$ as well. It is
not necessary that $\bar z$ correspond to parameters in the original
macroscopic range.

For $z\in K$ define 
\[
 \mathsf m(z)=\int_{\R^3}\phi(v)e^{z\cdot\phi(v)}  d v,
 \qquad
 \nabla_z \mathsf m(z)=\int_{\R^3}\phi(v)\phi(v)^{\mathsf T}
                                  e^{z\cdot\phi(v)}  d v.
\]
The Gaussian integrals and their derivatives are uniformly finite on
$K$. Also, using \eqref{eq:Gram-equivalence},
\[
 \xi^{\mathsf T}\nabla_z\mathsf m(z)\xi
 \ge\left(\inf_{K\times B}e^{z\cdot\phi(v)}\right)
                \lambda_B|\xi|^2.
\]
It follows that constants $\lambda_*,\Lambda_*>0$ satisfy
$\lambda_*I\le \nabla_z\mathsf m(z)\le\Lambda_*I$ on $K$. Since $K$ is convex, integration
along segments in $K$ gives
\begin{equation}\label{eq:moment-map}
 |\mathsf m(z)-\mathsf m(z')|\le\Lambda_*|z-z'|,
 \qquad
 |z-z'|\le\lambda_*^{-1}|\mathsf m(z)-\mathsf m(z')|.
\end{equation}
For the second bound, we took the scalar product of the difference of
moments with $z-z'$ before applying Cauchy--Schwarz to obtain 
the inverse estimate. 

Local moment matching and global conservation give
\[
 L^{-1}\int_{I\times\T^3}\mathsf m(z)  d x  d t
 =\mathsf m(z_\mu)=(1,0,3).
\]
Consequently,
\[
 |\mathsf m(\bar z)-\mathsf m(z_\mu)|
 \le\Lambda_*L^{-1/2}\norm{z-\bar z}_{L^2(I\times\T^3)}.
\]
Using \eqref{eq:moment-map} and the variance identity yields
\begin{align*}
 \norm{z-z_\mu}_2^2
 &=\norm{z-\bar z}_2^2+L|\bar z-z_\mu|^2\\
 &\le\bigl(1+(\Lambda_*/\lambda_*)^2\bigr)
                          \norm{z-\bar z}_2^2.
\end{align*}
This proves \eqref{eq:anchoring}.

For \eqref{eq:macro-equivalence}, we consider Gaussian integration to obtain
\begin{equation}\label{eq:macro-entropy}
	\begin{split}
 \HH_{\mac} &=\int_{\T^3}
 \left[\rho\log\rho-\rho+1 +\frac\rho2|U|^2
       +\frac{3\rho}{2}(T-1-\log T)\right]  d x  \\
       &:=\int_{\T^3}
       \left[ w(\rho) + \frac\rho2|U|^2
       +\frac{3\rho}{2} \psi(T) \right]  d x.
\end{split}
\end{equation}
We give quantitative bounds without assuming that $\rho$ or $T$
is close to $1$. Let $w(r)=r\log r-r+1$ and
$\psi(r)=r-1-\log r$. Both functions and their first derivatives
vanish at $1$, while $w^{\prime\prime}(r)=r^{-1}$ and $\psi''(r)=r^{-2}$.
Taylor's formula between $1$ and the indicated value yields
\begin{align*}
 \frac{|\rho-1|^2}{2\max\{1,\rho_+\}}
 &\le w(\rho)
 \le\frac{|\rho-1|^2}{2\min\{1,\rho_-\}},\\
 \frac{|T-1|^2}{2\max\{1,T_+\}^2}
 &\le\psi(T)
 \le\frac{|T-1|^2}{2\min\{1,T_-\}^2}.
\end{align*}
Together with $\rho_-\le\rho\le\rho_+$, these estimates bound
the integrand in \eqref{eq:macro-entropy} above and below by fixed
positive multiples of $|\rho-1|^2+|U|^2+|T-1|^2$. In particular,
one may take the lower constant to be
\[
 c_{\rm mac}=\min\left\{
 \frac1{2\max\{1,\rho_+\}},\frac{\rho_-}{2},
 \frac{3\rho_-}{4\max\{1,T_+\}^2}\right\}>0.
\]
The upper constant $C_{\rm mac}$ depends only on
$\rho_\pm,T_\pm$.

It remains to compare the physical variables with the coefficient
vector. The inverse formulas are
\[
 T=-\frac1{2c},\qquad U=-\frac b{2c}=Tb,\qquad
 \rho=\exp\!\left(a+\tfrac32\log(2\pi T)
                         +\frac{|U|^2}{2T}\right).
\]
The compact convex box $K$ has a fixed negative upper bound on $c$.
Thus the first two formulas, together with $z_\mu=(a_\mu,0,-1/2)$,
give
\begin{equation}\label{eq:app-coordinate-TU}
 |T-1|\le C|c+1/2|,\qquad |U|\le C|b|.
\end{equation}
Since $a_\mu=-\tfrac32\log(2\pi)$, taking logarithms in the
third formula gives
\[
 \log\rho=(a-a_\mu)+\tfrac32\log T+\frac{|U|^2}{2T}.
\]
The bounds $T_-\le T\le T_+$ imply
$|\log T|\le\min\{1,T_-\}^{-1}|T-1|$; the bounds on $U$ and $T$
give $|U|^2/(2T)\le (U_+/(2T_-))|U|$. Hence
\[
 |\log\rho|\le C\bigl(|a-a_\mu|+|b|+|c+1/2|\bigr)
              \le C|z-z_\mu|.
\]
Finally, the mean value theorem applied to the exponential, between
$0$ and $\log\rho$, shows
$|\rho-1|\le\max\{1,\rho_+\}|\log\rho|$. Therefore
\[
 |\rho-1|+|U|+|T-1|\le C|z-z_\mu|.
\]
Squaring, integrating, and combining with the preceding entropy
bounds proves \eqref{eq:macro-equivalence}.
\end{proof}

The following proposition provides finite-time control of the total
relative entropy under the source estimate
$\norm{\Tr g(t)}_{L^2(\T^3\times B)}^2
\le C_{\rm src}\DD(t)$.
The argument relies only on the transport structure and the
assumptions stated above, rather than on the specific form of the
Boltzmann collision operator. It therefore applies more generally
to kinetic transport equations satisfying these assumptions.
\begin{proposition}\label{prop:window}
Suppose the envelopes, conservation, and entropy decomposition hold.
Let ${\color{black}\DD(t)}\ge0$ be a measurable function, and assume that
$g=\log F$ satisfies
\begin{equation}\label{eq:abstract-source}
 \Tr g\in L^2_{\rm loc}([t_*,\infty)\times\T^3\times B),\qquad
 \norm{\Tr g(t)}_{L^2(\T^3\times B)}^2\le C_{\rm src}  \DD(t).
\end{equation}
Then, for $T_0$ from Lemma~\ref{lem:transport},
\begin{equation}\label{eq:abstract-window}
 \int_t^{t+T_0}\HH(s)  d s
 \le C_{\rm tran}\int_t^{t+T_0}
                         (\HH_{\mic}(s)+  \DD(s))  d s.
\end{equation}
The constant depends only on $C_{\rm tran}$, the envelopes, and the
geometric constants already specified.
\end{proposition}
\begin{proof}
Fix $I=(t,t+T_0)$ and write $P=z\cdot\phi=\log M_F$.
Because $P(s,x,\cdot)\in\cP$,
$\norm{(I-\Pi)g}_2\le\norm{g-P}_2$.
Lemmas~\ref{lem:transport}, ~\ref{lem:log-distance}, and \eqref{eq:abstract-source} give
\begin{align*}
 \inf_{P_0\in\cP}\norm{P-P_0}_2^2
 &\le 2\norm{P-g}_2^2+
                         2\inf_{P_0\in\cP}\norm{g-P_0}_2^2\\
 &\le 2\norm{P-g}_2^2+
  C\bigl(\norm{(I-\Pi)g}_2^2+T_0^2\norm{\Tr g}_2^2\bigr)\\
 &\le C_{T_0}\int_I\bigl(\HH_{\mic}(s)+
                         C_{\rm src}  \DD(s)\bigr)  d s,
\end{align*}
where all norms on the left are over $I\times\T^3\times B$.
Apply \eqref{eq:Gram-equivalence} and the fact that the best constant
$L^2$ approximation to $z$ is its space-time mean:
\[
 \norm{z-\bar z}_{L^2(I\times\T^3)}^2
 =\inf_{z_0\in\R^5}\norm{z-z_0}_{L^2(I\times\T^3)}^2
 \le C\int_I(\HH_{\mic}+  \DD)  d s.
\]
Lemma~\ref{lem:anchoring} now bounds $\int_I\HH_{\mac}$ by the same
right-hand side. Add $\int_I\HH_{\mic}$ and use
\eqref{eq:intro-split}.
\end{proof}

\begin{remark}
The velocity ball is used to identify the five coefficients of the
local Maxwellian. The microscopic entropy in
\eqref{eq:abstract-window} is always integrated over the whole
velocity space. Thus no tail error is discarded in passing from
logarithmic approximation on $B$ to global entropy.
\end{remark}

\section{The Boltzmann equation}\label{sec:boltzmann}

\subsection{Logarithmic source and entropy balance}
Write $Q=Q^+-Q^-$, with $Q^+,Q^-\ge0$ including the full factors of
$F$. Recall definition of $\Psi$ from \eqref{eq:Psi}. At fixed $(t,x)$ let
\[
 A=F'F_*',\quad Z=FF_*,\qquad
 \mathfrak d(v)=\int_{\R^3\times\Sph^2}
                    B(v-v_*,\sigma)\Psi(A,Z)  d\sigma  d v_*.
\]
Then $\int\mathfrak d(v)  d v=4D_v(F)$. 

\begin{lemma} 
	\label{lem:BE-source}
Under \eqref{eq:kernel} and \eqref{eq:envelope}, for every fixed
$R<\infty$ there is $K_R<\infty$ such that
\begin{equation}\label{eq:Qbound}
 Q^+(F,F)(v)+Q^-(F,F)(v)\le K_R,\qquad |v|\le R.
\end{equation}
Moreover,
\begin{equation}\label{eq:Q-square}
 |Q(F,F)(v)|^2\le\tfrac12(Q^++Q^-)(v)\mathfrak d(v),
\end{equation}
and, with $\ell_R=a_*e^{-b_*R^2}$,
\begin{equation}\label{eq:QoverF}
 \int_{B_R}|Q(F,F)/F|^2  d v
 \le\frac{2K_R}{\ell_R^2}D_v(F).
\end{equation}
\end{lemma}
\begin{proof}
Energy conservation implies
\[
 \av{v'}^2\av{v_*'}^2\ge1+|v'|^2+|v_*'|^2
 =1+|v|^2+|v_*|^2\ge\av{v_*}^2.
\]
Hence $F'F_*'\le M_0^2\av{v_*}^{-m}$, and the same upper bound
holds for $FF_*$ when $v$ is arbitrary. Also
$|v-v_*|^\gamma\le2^{\gamma/2}\av v^\gamma\av{v_*}^\gamma$.
One may therefore take
\begin{equation}\label{eq:KR}
 K_R=2^{1+\gamma/2}\beta_bM_0^2\av R^\gamma
                         \int_{\R^3}\av{v_*}^{\gamma-m}  d v_*.
\end{equation}
This requires only $m>3+\gamma$.

For $A,Z>0$ the elementary inequality
\[
 (A-Z)\log(A/Z)\ge\frac{2(A-Z)^2}{A+Z}
\]
follows by differentiating $\log r-2(r-1)/(r+1)$ for $r\ge1$ and
then interchanging $A,Z$ if needed. Weighted Cauchy--Schwarz gives
\begin{align*}
 |Q(v)|^2
 &\le\left(\int B (A+Z)\right)
             \left(\int B \frac{(A-Z)^2}{A+Z}\right)\\
 &\le\tfrac12(Q^++Q^-)(v)\mathfrak d(v),
\end{align*}
where $B=B(v-v_*,\sigma)$ is collision kernel.
All omitted integrals here are in $(v_*,\sigma)$. The first two
factors are finite by \eqref{eq:Qbound} and
$(A-Z)^2/(A+Z)\le A+Z$. Finally divide by $F^2\ge\ell_R^2$ and
integrate on $B_R$, using nonnegativity of $\mathfrak d$ and its
full-velocity integral $4D_v(F)$.
\end{proof}

\begin{lemma}\label{lem:BE-chain}
For the mild solutions of Theorem~\ref{thm:boltzmann},
\begin{equation}\label{eq:BE-logchain}
 \Tr\log F=Q(F,F)/F
\end{equation}
in distributions on every finite cylinder with bounded velocities.
Furthermore, $\HH$ is locally absolutely continuous and
\begin{equation}\label{eq:BE-entropy}
 \HH(t)+\int_s^t  \DD_{\rm B}(\tau)  d\tau=\HH(s),
 \qquad t\ge s\ge t_*.
\end{equation}
\end{lemma}
\begin{proof}
Fix $I=(s_0,s_1)$ and $R<\infty$, and put
$F^\#(t,y,v)=F(t,y+(t-s_0)v,v)$.
The mild formula and \eqref{eq:Qbound} give
$\partial_tF^\#=Q(F,F)^\#$ with bounded right-hand side on
$I\times\T^3\times B_R$. In particular $F^\#\in H^1_tL^2_{y,v}$,
and for almost every $(y,v)$ it has an absolutely continuous
representative. Its almost-everywhere lower bound $\ell_R$ on $B_R$ extends
to every time for that representative. The one-dimensional chain
rule therefore gives
\[
 \log F^\#(t)-\log F^\#(s)
 =\int_s^t\frac{Q(F,F)^\#}{F^\#}(\tau)  d\tau.
\]
Changing coordinates back proves \eqref{eq:BE-logchain}.

Now we prove \eqref{eq:BE-entropy}. For the entropy identity, \eqref{eq:log-bound} gives
$|\log(F/\mu)|\le C\av v^2$. By the collision symmetries and
energy conservation,
\begin{equation}\label{eq:weighted-Q}
 \int\av v^2(Q^++Q^-)  d v
 \le C\left(\int\av v^{2+\gamma}F  d v\right)
       \left(\int\av v^\gamma F  d v\right),\quad m>5+\gamma.
\end{equation}
For instance, transfer the gain term to pre-collisional velocities,
bound $\av{v'}^2$ by $\av v^2+\av{v_*}^2$, and use the same
bound for $|v-v_*|^\gamma$ as above. The two terms are equal after
interchanging $v,v_*$. The right side is uniform because $m>5+\gamma$.
Thus $Q\log(F/\mu)$ is integrable on every finite space-time cylinder
with the full velocity space. Apply the characteristic chain rule to
$\Phi(F\mid\mu)$ on $B_R$, integrate in $x,v$, and let $R\to\infty$.
Dominated convergence for the source and integrability of the entropy
density justify the limit. Finally the collision-invariant part
$\log\mu$ integrates to zero, and symmetrization gives
\[
 \int_{\R^3}Q(F,F)\log F  d v=-D_v(F).
\]
These manipulations are absolutely integrable: each logarithmic
factor at a pre- or post-collisional velocity is controlled by
$C(1+|v|^2+|v_*|^2)$, and \eqref{eq:weighted-Q} applies.
The resulting integral identity is \eqref{eq:BE-entropy}.
\end{proof}

The proof of decay below would also apply if
\eqref{eq:BE-entropy} were replaced by the entropy inequality.
The stronger weight $m>5+\gamma$ is used here to justify the identity
directly from the mild formulation. The logarithmic source bound on
the fixed ball only needs the weaker local collision integrability.

\begin{lemma} 
\label{lem:weighted-to-entropy}
If $m>5$, $0\le f(x,v)\le M\av v^{-m}$, and
$N=\norm{\av v^m(f-\mu)}_{L^\infty_{x,v}}<\infty$, then
\begin{equation}\label{eq:weighted-to-entropy}
 H(f\mid\mu)\le C_{m,M}N.
\end{equation}
The spatial torus has normalized measure.
\end{lemma}
\begin{proof}
On $f<\mu$, one has $\Phi(f\mid\mu)\le\mu-f$. On $f\ge\mu$,
$r-1-\log r\ge0$ for $r=f/\mu\ge1$ gives
$\Phi(f\mid\mu)\le(f-\mu)\log(f/\mu)$.
The upper envelope implies
$0\le\log(f/\mu)\le C_M(1+|v|^2)$ on this latter set.
Thus, in both cases,
\[
 0\le\Phi(f\mid\mu)\le C_M(1+|v|^2)|f-\mu|.
\]
After integration, \eqref{eq:weighted-to-entropy} follows from
$\int_{\R^3}(1+|v|^2)\av v^{-m}  d v<\infty$ for $m>5$.
\end{proof}

\subsection{Proofs of the Boltzmann decay statements}

\begin{proof}[Proof of Theorem~\ref{thm:boltzmann}]
We first establish the finite-time estimate.
Integrate \eqref{eq:QoverF} with $R=1$ in $x$ and use
\eqref{eq:BE-logchain}:
\begin{equation}\label{eq:BE-source-time}
 \norm{\Tr\log F(t)}_{L^2(\T^3\times B)}^2
 \le\frac{2K_1}{\ell_1^2}  \DD_{\rm B}(t).
\end{equation}
By \eqref{eq:BE-entropy}, the right-hand side is integrable on
every finite time interval. Thus all hypotheses of
Proposition~\ref{prop:window} hold, and it gives
\eqref{eq:BE-window}. This proves the first assertion of
Theorem~\ref{thm:boltzmann}, without using \eqref{eq:MI}.

We next derive the algebraic decay.
Assume the additional uniform angular lower bound $b\ge c_b>0$
almost everywhere and $m>7$ in Theorem~\ref{thm:boltzmann}. Proposition~\ref{aux:proposition}
then proves \eqref{eq:MI} uniformly with $p=1+\varepsilon_*$.
Let
$\Delta_t=\int_t^{t+T_0} {\color{black}\DD_{\rm B}(s)}  d s$. H\"older's
inequality on a cylinder of measure $T_0$ gives
\begin{equation}\label{eq:mic-time}
 \int_t^{t+T_0}\HH_{\mic}(s)  d s
 \le T_0^{1-1/p}c_{\mic}^{-1/p}\Delta_t^{1/p}.
\end{equation}
Monotonicity of entropy and \eqref{eq:BE-window} imply
\[
 \HH(t+T_0)\le A_0\Delta_t^{1/p}+A_1\Delta_t,
 \quad A_0=C_{\rm B}T_0^{-1/p}c_{\mic}^{-1/p},\quad
 A_1=C_{\rm B}/T_0.
\]
Let $H_*=\HH(t_*)$. If $H_*=0$, then $F=\mu$ at all subsequent
times. Otherwise $\Delta_t\le H_*$, and consequently
\begin{equation}\label{eq:finite-decay}
 \HH(t+T_0)^p\le C_*\Delta_t
 \le C_*\bigl(\HH(t)-\HH(t+T_0)\bigr),
 \quad C_*=(A_0+A_1H_*^{1-1/p})^p.
\end{equation}
This is the finite-time entropy--dissipation inequality.

To solve the recurrence with the future-time power in
\eqref{eq:finite-decay}, let
$y_n=\HH(t_*+nT_0)$, $c_*=C_*^{-1}$, and $r=1/(p-1)$. Then
\begin{equation}\label{eq:recurrence}
 y_n\ge y_{n+1}+c_*y_{n+1}^p.
\end{equation}
Choose
\[
 A\ge H_*,\qquad A^{p-1}\ge r2^{r+1}/c_*,\qquad
 b_n=A(n+1)^{-r}.
\]
The mean value theorem gives
\[
 b_n-b_{n+1}\le Ar(n+1)^{-r-1}
 \le Ar2^{r+1}(n+2)^{-r-1}\le c_*b_{n+1}^p.
\]
Since $u\mapsto u+c_*u^p$ is increasing, \eqref{eq:recurrence}
and $y_0\le b_0$ imply $y_n\le b_n$ by induction. For
$n=\lfloor(t-t_*)/T_0\rfloor$, monotonicity gives
\begin{equation}\label{eq:continuous-polynomial}
 \HH(t)\le A(n+1)^{-r}
 \le A(1+T_0)^r(1+t-t_*)^{-r}.
\end{equation}
Kullback-Pinsker's inequality, with both total masses equal to one, gives
$\norm{F-\mu}_{L^1}^2\le2\HH$.
\end{proof}

\begin{proof}[Proof of Corollary~\ref{coro:boltzmann}]
The extra assumptions imply $m>7$ and $b\ge {\color{black}c_b}$, so
Theorem~\ref{thm:boltzmann} gives $\HH(t)\to0$ at an explicit
algebraic rate. Put $w(v)=\av v^m$. The uniform upper envelope gives
\[
 \sup_{t\ge t_*}\norm{w(F(t)-\mu)}_{L^\infty_{x,v}}
 \le M_0+\norm{w\mu}_{L^\infty_v}<
 \overline M:=1+M_0+\norm{w\mu}_{L^\infty_v}.
\]
Use $M_1=\overline M\ge1$. The polynomial large-amplitude
theorem \cite[Theorem 7.26]{HJJLP26}, with weight exponent
$m>\beta_\infty$, supplies a positive small-entropy threshold
$\varepsilon_0$ for the fixed bounds $M_1,\overline M$.
Choose $s\ge t_*$ so that $\HH(s)\le\varepsilon_0$; the
algebraic estimate makes such an $s$ finite. The translated solution
$F(s+\cdot)$ has the same conserved moments, the required initial
weighted norm bound, and the global future bound $\overline M$.
The Theorem 7.26 of \cite{HJJLP26} therefore yields
\[
 \norm{w(F(t)-\mu)}_{L^\infty_{x,v}}
 \le C_A e^{-\lambda_\infty(t-s)},\qquad t\ge s.
\]
On $t_*\le t\le s$ the bound $\overline M$ holds. Enlarging the
prefactor by a finite factor depending on $s-t_*$ gives the first
estimate in \eqref{eq:BE-exp-decay} for every $t\ge t_*$.
Finally Lemma~\ref{lem:weighted-to-entropy}, applied to the
uniform upper envelope of $F(t)$, gives
\[
 \HH(t)\le C_{m,M_0}\norm{w(F(t)-\mu)}_{L^\infty_{x,v}}
 \le C_H e^{-\lambda_\infty(t-t_*)}.
\]
This conversion uses only $m>5$ and does not require a separate
smallness assumption on the weighted norm.
\end{proof}

\begin{corollary} 
	\label{cor:BE-macro}
Under the additional assumptions $b\ge c_b>0$ almost everywhere
and $m>7$ in Theorem~\ref{thm:boltzmann},
\[
 \norm{\rho(t)-1}_{L^2_x}^2+\norm{U(t)}_{L^2_x}^2
 +\norm{T(t)-1}_{L^2_x}^2\le C(1+t-t_*)^{-1/(p-1)}
\]
at times for which the macroscopic fields are defined.
\end{corollary}
\begin{proof}
Use the lower bound in \eqref{eq:macro-equivalence} and
$\HH_{\mac}\le\HH$.
\end{proof}

\begin{remark}[Uniformity and scope of microscopic coercivity]
The constant $c_{\mic}$ in \eqref{eq:MI} is uniform over the
time-space slices because all of them share the same envelopes.
The algebraic estimate alone does not yield exponential entropy
decay; Corollary~\ref{coro:boltzmann} additionally invokes the
large-amplitude weighted theorem. A fixed polynomial exponent $m$
supplies moments only of orders less than $m-3$, so it does not
provide arbitrary microscopic exponents by increasing moment order. \\
\end{remark}

\section{The BGK equation: exponential entropy decay}
\label{sec:bgk}

We turn to the normalized BGK equation
\begin{equation}\label{bgk:equation}
 \partial_tF+v\cdot\nabla_xF=M_F-F.
\end{equation}
The collision frequency is one. Because the BGK entropy production
controls microscopic entropy linearly, the finite-time transport
estimate will give exponential relaxation.

\subsection{Entropy production and the logarithmic equation}

Define
\begin{equation}\label{bgk:dissipation}
 \mathcal D_{\rm G}(t)
 :=\int_{\mathbb T^3\times\mathbb R^3}
 (F-M_F)\log\frac{F}{M_F}\,dx\,dv.
\end{equation}
We retain the notation
$\mathcal H=H(F\mid\mu)$,
$\mathcal H_{\rm mic}=H(F\mid M_F)$, and
$\mathcal H_{\rm mac}=H(M_F\mid\mu)$.

\begin{lemma}\label{bgk:entropy}
Let $F$ be a positive mild solution of \eqref{bgk:equation} satisfying \eqref{eq:envelope} with $m>5$.
Then $\mathcal H$ is locally absolutely continuous and
\begin{equation}\label{bgk:entropy-identity}
 \mathcal H(t)+\int_s^t\mathcal D_{\rm G}(\tau)\,d\tau
 =\mathcal H(s),\qquad t_*\le s\le t.
\end{equation}
Moreover, for almost every $t$,
\begin{equation}\label{bgk:linear-micro}
 \mathcal D_{\rm G}(t)
 =H(F(t)\mid M_F(t))+H(M_F(t)\mid F(t))
 \ge\mathcal H_{\rm mic}(t).
\end{equation}
On the fixed velocity ball $B=B_1(0)$, the function $g=\log F$ satisfies
\begin{equation}\label{bgk:log-equation}
 \mathcal T g=\frac{M_F-F}{F}
 \quad\hbox{in distributions},\qquad
 \mathcal T=\partial_t+v\cdot\nabla_x.
\end{equation}
\end{lemma}

\begin{proof}
The uniform macroscopic bounds give constants $C,c>0$ such that
$M_F(t,x,v)\le C e^{-c|v|^2}$.  They also imply
$|\log M_F(t,x,v)|\le C(1+|v|^2)$.  The Gaussian lower bound for $F$
and its uniform upper bound imply
$|\log F(t,x,v)|\le C(1+|v|^2)$.
Consequently, 
\[
\begin{aligned}
	\left|(M_F-F)\log\frac{F}{\mu}\right|
	&\le (M_F+F)\left|\log\frac{F}{\mu}\right|\\
	&\le C_2(1+|v|^2)
	\left(C_1e^{-c|v|^2}+M_0\langle v\rangle^{-m}\right)\\
	&\le C(1+|v|^2)
	\left(M_0\langle v\rangle^{-m}+e^{-c|v|^2}\right),
\end{aligned}
\]
and is integrable since $m>5$.

Along almost every free characteristic,
$f(t)=F(t,x+tv,v)$ is locally absolutely continuous and satisfies
\[
 f'(t)=M_F(t,x+tv,v)-f(t).
\]
For each fixed $v$, the lower bound for $F$ stays strictly positive.
The scalar chain rule therefore applies both to $\log f$ and to
$f\log(f/\mu)-f+\mu$.  The first gives \eqref{bgk:log-equation}.
For the second, integrate along characteristics and then in $(x,v)$.
Translation invariance of the torus and the preceding domination
justify Fubini's theorem and give
\[
 \mathcal H(t)-\mathcal H(s)
 =\int_s^t\int_{\mathbb T^3\times\mathbb R^3}
 (M_F-F)\log\frac{F}{\mu}\,dx\,dv\,d\tau.
\]
Since $\log(M_F/\mu)$ is a linear combination of
$1,v_1,v_2,v_3,|v|^2$, local moment matching yields
\[
 \int_{\mathbb R^3}(M_F-F)\log\frac{M_F}{\mu}\,dv=0
 \quad\hbox{for almost every }(t,x).
\]
Subtracting this identity proves \eqref{bgk:entropy-identity}.
Finally, with the convention
$\Phi(a\mid b)=a\log(a/b)-a+b$, the pointwise identity
\[
 \Phi(a\mid b)+\Phi(b\mid a)=(a-b)(\log a-\log b)
\]
and the nonnegativity of both relative entropies prove
\eqref{bgk:linear-micro}.
\end{proof}

\begin{lemma}\label{bgk:log-control}
Let $m_B,M_B$ be positive constants such that
\begin{equation}\label{bgk:ball-bounds}
 0<m_B\le F(t,x,v),M_F(t,x,v)\le M_B,
 \qquad v\in B.
\end{equation}
With $P=\log M_F$ and with $\Pi$ the orthogonal projection onto
$\mathcal P=\operatorname{span}\{1,v_1,v_2,v_3,|v|^2\}$ in $L^2(B)$,
one has
\begin{align}
 \|g-P\|_{L^2(\mathbb T^3\times B)}^2
 &\le \frac{1}{m_B}\mathcal D_{\rm G}(t),
 \label{bgk:log-distance}\\
 \|(I-\Pi)g\|_{L^2(\mathbb T^3\times B)}^2
 &\le \frac{1}{m_B}\mathcal D_{\rm G}(t),
 \label{bgk:projection}\\
 \|\mathcal T g\|_{L^2(\mathbb T^3\times B)}^2
 &\le \frac{M_B}{m_B^2}\mathcal D_{\rm G}(t).
 \label{bgk:source}
\end{align}
\end{lemma}

\begin{proof}
For $a,b\in[m_B,M_B]$, the mean value theorem gives
$\log a-\log b=(a-b)/\xi$ for some
$\xi\in[m_B,M_B]$.  Thus
\begin{align*}
 |\log a-\log b|^2
 &\le\frac{1}{m_B}(a-b)(\log a-\log b),\\
 |a-b|^2
 &\le M_B(a-b)(\log a-\log b).
\end{align*}
The first inequality, with $a=F$ and $b=M_F$, proves
\eqref{bgk:log-distance}.  Since $P(t,x,\cdot)\in\mathcal P$,
\[
 (I-\Pi)g=(I-\Pi)(g-P),
\]
and the contraction property of an orthogonal projection gives
\eqref{bgk:projection}.  Finally, \eqref{bgk:log-equation} and the
second scalar inequality imply
\[
 \left|\frac{M_F-F}{F}\right|^2
 \le \frac{M_B}{m_B^2}
 (F-M_F)\log\frac{F}{M_F}.
\]
Integration over $\mathbb T^3\times B$ proves \eqref{bgk:source},
because the dissipation integrand is nonnegative on the whole velocity
space.
\end{proof}

\subsection{Proofs of the BGK decay statements}

\begin{proof}[Proof of Theorem~\ref{thm:bgk}]
The uniform pointwise bounds provide the macroscopic bounds and
\eqref{bgk:ball-bounds}.  Lemma~\ref{bgk:log-control} verifies the
logarithmic transport hypotheses in the common argument.  Applying
Proposition~\ref{prop:window}, followed by
\eqref{bgk:linear-micro}, gives
\[
 \int_t^{t+T_0}\mathcal H(s)\,ds
 \le C\int_t^{t+T_0}
 \bigl(\mathcal H_{\rm mic}(s)+\mathcal D_{\rm G}(s)\bigr)\,ds
 \le K_{\rm G}\int_t^{t+T_0}\mathcal D_{\rm G}(s)\,ds.
\]
Equivalently, one can apply Lemma~\ref{lem:transport} to $g$ using
\eqref{bgk:projection} and \eqref{bgk:source}, replace $g$ by $P$
using \eqref{bgk:log-distance}, and then use conservation of the
global moments through Lemma~\ref{lem:anchoring}.  This directly
controls the time integral of $\mathcal H_{\rm mac}$ by the same
dissipation integral.  Thus both microscopic and macroscopic parts
are controlled linearly in the BGK case.

By monotonicity of $\mathcal H$ and \eqref{bgk:entropy-identity},
\begin{align*}
 T_0\mathcal H(t+T_0)
 &\le\int_t^{t+T_0}\mathcal H(s)\,ds\\
 &\le K_{\rm G}\bigl(\mathcal H(t)-\mathcal H(t+T_0)\bigr).
\end{align*}
With $q=K_{\rm G}/(K_{\rm G}+T_0)$ and
$\lambda=T_0^{-1}\log(1/q)$, it follows that
\begin{equation}\label{bgk:contraction}
 \mathcal H(t+T_0)\le q\mathcal H(t).
\end{equation}
For $n=\lfloor(t-t_*)/T_0\rfloor$, iteration and monotonicity imply
\[
 \mathcal H(t)\le q^n\mathcal H(t_*)
 \le q^{-1}e^{-\lambda(t-t_*)}\mathcal H(t_*).
\]
Since the entropy splitting has nonnegative terms,
\[
 0\le\mathcal H_{\rm mic}(t),\mathcal H_{\rm mac}(t)
 \le\mathcal H(t),
\]
both microscopic and macroscopic entropies have the same exponential
upper bound, without requiring either term to be monotone.
Kullback-Pinsker's inequality gives the stated $L^1$ estimate. 

We next prove \eqref{eq:BGK-decay-norm}. Assume $m>10$ and put $w(v)=\av v^m$. The envelope implies
\[
 \sup_{t\ge t_*}\norm{w(F(t)-\mu)}_{L^\infty_{x,v}}
 \le M_0+\norm{w\mu}_{L^\infty_v}<
 \overline M:=1+M_0+\norm{w\mu}_{L^\infty_v}.
\]
For every restart time $s\ge t_*$, the Gaussian lower bound gives
the free-streaming density condition in \cite[Theorem 1.1]{BKLY26}:
\[
 \inf_{\tau\ge0,\;x\in\T^3}
 \int_{\R^3}F(s,x-\tau v,v)  d v
 \ge a_*\int_{\R^3}e^{-b_*|v|^2}  d v=:c_{\rm fs}>0.
\]
The weight $1+|v|^m$ used in that paper is equivalent to our
$\av v^m$. Its collision-frequency class $\rho^aT^b$ contains
the unit-frequency equation by taking $a=b=0$, and its required
weight range is $m>10$.

The entropy bound \eqref{eq:BGK-decay-main} permits us to choose a
finite $s\ge t_*$ for which $\HH(s)$ is below the small-entropy
threshold of \cite[Theorem 1.1]{BKLY26} corresponding to
$\overline M$ and $c_{\rm fs}$. Conservation supplies the normalized
equilibrium, and the translated mild solution belongs to the
uniqueness class of the cited result. It follows that
$\norm{w(F(t)-\mu)}_\infty\le C_s e^{-\kappa(t-s)}$ for $t\ge s$.
The uniform bound $\overline M$ on $[t_*,s]$ absorbs the finite
restart shift into the constant and proves
\eqref{eq:BGK-decay-norm}. Lemma~\ref{lem:weighted-to-entropy}
then gives $\HH(t)\le C_{m,M_0}C_\infty
e^{-\kappa(t-t_*)}$, as claimed.
\end{proof}

\begin{remark}\label{bgk:scope}
The BGK entropy estimate assumes two uniform velocity envelopes
for a global mild solution. No spatial or velocity Sobolev regularity
enters its proof.
In fact, the upper-bound hypothesis in this section can be replaced by
\[
F(t,x,v)\le W(v),\qquad
W\in L^\infty(\mathbb R^3),\qquad
\int_{\mathbb R^3}(1+|v|^2)W(v)\,dv<\infty,
\]
while retaining the
Gaussian lower bound.  The conclusions concern entropy, $L^1$, and
the macroscopic $L^2$ fields. The weighted $L^\infty_{x,v}$ conclusion requires, in addition,
the polynomial upper envelope with $m>10$ and the
large-amplitude result of \cite{BKLY26}. \\
\end{remark}

\section{Application to the Perthame--Pulvirenti BGK solution class}
\label{sec:application}

Perthame and Pulvirenti \cite{PP93} constructed a global unique
BGK mild solution in a polynomially weighted $L^\infty$ space.
Their estimate holds on every finite time interval, with a constant
that may depend on the interval. Here we assume in addition that
the local density stays bounded above. We show that this one
uniform bound leads to the envelopes of Theorem~\ref{thm:bgk},
without a small-entropy assumption.

\begin{proof} [Proof of Theorem \ref{thm:PP-application}]
The lower bound in \eqref{eq:app-data} gives the free-transport density
condition of \cite[Theorem 3.1]{PP93}:
\[
 \int_{\R^3}F_0(x-\tau v,v)  d v
 \ge\int_{\R^3}\varphi(v)  d v>0
 \quad\text{for every }(x,\tau).
\]
That theorem provides the stated global mild solution and a finite
weighted $L^\infty$ bound on every finite time interval. Conservation
holds. Moreover, its finite-horizon bounds justify the BGK $H$-theorem
first away from $t=0$, and then by passage to $t=0$:
\begin{equation}\label{eq:app-entropy-monotone}
 0\le\HH(t)\le H_0:=\HH(0)<\infty.
\end{equation}
One may use the entropy inequality here; an entropy identity is not
needed in the present step.

We begin with lower bounds on $\rho$ and $T$ after a fixed
waiting period. All spatial measures
in this argument are normalized so that $|\T^3|=1$. Write
\[
 e_2(s,x)=\int_{\R^3}|v|^2F(s,x,v)  d v
          =\rho(s,x)\bigl(|U(s,x)|^2+3T(s,x)\bigr).
\]
Conservation gives $\int\rho  d x=1$ and $\int e_2  d x=3$.
Set $A_{\rho} \ge1$ and $g_0=(2A_{\rho})^{-1}$. If we define 
$E_s=\{x:\rho(s,x)\ge1/2\}$, then
\[
 1=\int_{\T^3}\rho  d x
 \le\tfrac12(1-|E_s|)+A_{\rho}|E_s|
 \le\tfrac12+A_{\rho}|E_s|, \ \ \text{hence} \ \ |E_s|\ge g_0.
\]
Choose $E_*=24A_{\rho}$. Markov's inequality and energy conservation give
$|\{e_2(s,\cdot)>E_*\}|\le3/E_*=g_0/4$.
{\color{black}\hypersetup{linkcolor=black,citecolor=black}
Define $\psi(r)=r-1-\log r$ for $r>0$. This function is nonnegative,
decreases on $(0,1)$, and tends to infinity as $r\downarrow0$.
By the entropy splitting and \eqref{eq:macro-entropy},
\[
 \frac32\int_{\T^3}\rho(s,x)\psi(T(s,x))\,d x
 \le\HH_{\mac}(s)\le\HH(s)\le H_0.
\]
For $0<\theta<1$, on $E_s\cap\{x:T(s,x)<\theta\}$,
$\rho\ge1/2$ and $\psi(T)\ge\psi(\theta)$. Consequently,
\begin{align*}
 H_0
 &\ge\frac32\int_{\T^3}\rho(s,x)\psi(T(s,x))\,d x\\
 &\ge\frac32\psi(\theta)\int_{E_s\cap\{x:T(s,x)<\theta\}}\rho(s,x)\,d x\\
 &\ge\frac34\psi(\theta)|E_s\cap\{x:T(s,x)<\theta\}|,
\end{align*}
and therefore
\[
 |E_s\cap\{T<\theta\}|
 \le\frac{4H_0}{3\psi(\theta)}.
\]
}
Fix $\theta\in(0,1)$ so that
$\psi(\theta)\ge 16H_0/(3g_0)$ and hence $|E_s\cap\{T<\theta\}| \leq g_0/4$ (with any $\theta$ if $H_0=0$).
The good spatial set
\[
 G_s :=\{x:\rho(s,x)\ge 1/2,\ e_2(s,x)\le E_*,\
                         T(s,x)\ge\theta\}
\]
therefore satisfies $|G_s|\ge g_0 - \frac{g_0}{4} - \frac{g_0}{4} =g_0/2$ for every $s\ge0$.
No regularity or common location of the sets $G_s$ is used.

On $G_s$, $\rho\in[1/2,A_\rho]$ and
$|U|^2+3T=e_2/\rho\le2E_*$. Put
$U_*=\sqrt{2E_*}$ and $T_*=2E_*/3$. Then
$|U|\le U_*$ and $\theta\le T\le T_*$. Since
$|v-U|^2\le2|v|^2+2U_*^2$, the local Maxwellian obeys
\begin{equation}\label{eq:app-good-M}
 M_F(s,y,v)\ge c_0e^{-C_0|v|^2}\mathbf1_{G_s}(y),
 \quad c_0=\frac{e^{-U_*^2/\theta}}{2(2\pi T_*)^{3/2}},
 \quad C_0=\theta^{-1}.
\end{equation}
For $t\ge2$, keep only delays $\tau\in[1,2]$ in the positive
BGK Duhamel formula. This gives the pointwise lower function
\begin{equation}\label{eq:app-bath}
 B_t(x,v):=c_0e^{-C_0|v|^2}
       \int_1^2e^{-\tau}\mathbf1_{G_{t-\tau}}(x-\tau v)
         d\tau\le F(t,x,v).
\end{equation}
The next calculation shows that the mass and centered second moment
of $B_t$ are bounded below independently of $t$ and $x$.

Represent $\T^3$ by $[0,2\pi)^3$ and fix representatives
$x,y\in[0,2\pi)^3$. For $\tau\in[1,2]$ and
$n\in\mathcal N:=\{0,e_1,e_2,e_3\}$, set
\[
 v_n(y)=\frac{x-y+2\pi n}{\tau}.
\]
The map $y\mapsto v_n(y)$ has Jacobian $\tau^{-3}$, and
$x-\tau v_n(y)=y$ modulo $2\pi\mathbb Z^3$. Its four image cells are
pairwise disjoint up to null boundaries; all their velocities have
$|v_n(y)|\le R_*:=4\pi\sqrt3$. With normalized torus measure,
the corresponding change of variables is
$  d v=(2\pi)^3\tau^{-3}  d y$. Restricting the velocity integral
to these cells and using \eqref{eq:app-bath} therefore gives
\begin{equation} \label{eq:B lower bound}
	\begin{split}
 \int_{\R^3} B_t(x,v)  d v
 &\ge c_0e^{-C_0R_*^2}(2\pi)^3
   \int_1^2e^{-\tau}\tau^{-3}|G_{t-\tau}|  d\tau\\
 &\ge c_0e^{-C_0R_*^2}(2\pi)^3 \frac{g_0}{2}
   \int_1^2e^{-\tau}\tau^{-3}  d\tau>0,
\end{split}
\end{equation}
since $|G_s|\geq \frac{g_0}{2}$.

For the centered moment, let
$\bar v(y)=\frac14\sum_{n\in\mathcal N}v_n(y)$. The variance
identity and the geometry of the four lifts give, for any $a\in\R^3$,
\begin{align*}
 \sum_{n\in\mathcal N}|v_n(y)-a|^2
 &=4|a-\bar v(y)|^2
   +\sum_{n\in\mathcal N}|v_n(y)-\bar v(y)|^2\\
 &=4|a-\bar v(y)|^2+\frac{9\pi^2}{\tau^2}.
\end{align*}
Because $|\bar v(y)|\le R_*$ and $1\le\tau\le2$,
there is a numerical $c_{\rm geom}>0$, depending only on the
chosen torus, such that
$\sum_{n\in\mathcal N}|v_n(y)-a|^2
 \ge c_{\rm geom}(1+|a|^2)$.
Indeed $|a|^2\le2|a-\bar v|^2+2R_*^2$, while
$9\pi^2/\tau^2\ge9\pi^2/4$. The same change of variables yields
\begin{equation} \label{eq:B 2nd}
	\begin{split}
 \int_{\R^3}|v-a|^2B_t(x,v)  d v
 \ge c_0e^{-C_0R_*^2}(2\pi)^3c_{\rm geom}{\color{black}\frac{g_0}{2}}
    \left(\int_1^2e^{-\tau}\tau^{-3}  d\tau\right)(1+|a|^2).
\end{split}
\end{equation}
Choose $b_0>0$ no larger than the two displayed positive constants in \eqref{eq:B lower bound} and \eqref{eq:B 2nd}. Then we obtain
\begin{equation}\label{eq:app-bath-moments}
 \int B_t(x,v)  d v\ge b_0,\qquad
 \int|v-a|^2B_t(x,v)  d v\ge b_0(1+|a|^2)
 \quad(a\in\R^3).
\end{equation}
Since $B_t\le F$, taking $a=U(t,x)$ gives, for $t\ge2$,
\begin{equation}\label{eq:app-macro-lower}
 \rho(t,x)\ge b_0,\qquad
 3\rho(t,x)T(t,x)\ge b_0(1+|U(t,x)|^2),
 \qquad T(t,x)\ge\frac{b_0}{3A_\rho}.
\end{equation}
This step uses only the density bound, conservation, and entropy
monotonicity. It also shows that the Maxwellian height
$q_F=\rho(2\pi T)^{-3/2}$ is uniformly bounded.

We next show that the density bound makes the weighted estimate
of \cite{PP93} uniform in time. For a nonnegative velocity density
$h$, set
\[
 N_\beta(h)={\color{black}\sup_{v\in\R^3}}\av v^\beta h(v),\qquad
 P_\beta(h)=N_\beta(M_h).
\]
{\color{black}\hypersetup{linkcolor=black,citecolor=black}
When $h=h(x,v)$ also depends on position, we use the notation
\[
 N_\beta(h):=\sup_{x\in\T^3}
             N_\beta(h(x,\cdot)),\qquad
 P_\beta(h):=\sup_{x\in\T^3}
             P_\beta(h(x,\cdot)).
\]
}
We recall the weighted Maxwellian estimate needed below. Write
$N=N_\beta(h)$. If $h=0$,
set $M_h=0$ and there is nothing to prove. Otherwise let $\mathfrak{m}$, $u$,
and $\vartheta>0$ be its density, mean velocity, and temperature.
The estimate of Perthame--Pulvirenti \cite[Proposition 2.1 and
(2.11)]{PP93} is
\begin{equation}\label{eq:app-PP-weighted}
 P_\beta(h)=\sup_v\av v^\beta M_h(v)
 \le C_\beta N_\beta(h),\qquad \beta>5.
\end{equation}
We briefly explain \eqref{eq:app-PP-weighted} to see why the constant is independent of
$\mathfrak{m},u,\vartheta$. Since
$\int|v-u|^2h  d v=3\mathfrak{m}\vartheta$, at least half the mass of $h$
lies in the ball $|v-u|\le\sqrt{6\vartheta}$. If
$|u|\ge2\sqrt{6\vartheta}$, then $\av v\ge c\av u$ on that ball,
and its volume is $C\vartheta^{3/2}$. Thus
\[
 \mathfrak{m}\le C_\beta N\vartheta^{3/2}\av u^{-\beta}
 \quad\text{when }|u|\ge2\sqrt{6\vartheta}.
\]
If $|u|<2\sqrt{6\vartheta}$ and $\vartheta\le1$, the simple bound
$h\le N$ on the same ball gives $\mathfrak{m}\le C N\vartheta^{3/2}$.
Finally, if $|u|<2\sqrt{6\vartheta}$ and $\vartheta>1$, put
$R^2=3\vartheta/2$. Splitting the raw second moment at $|v|=R$
and using $\beta>5$ yields
\begin{align*}
 3\mathfrak{m}\vartheta
 &\le\int|v|^2h(v)  d v\\
 &\le R^2 \mathfrak{m}+N\int_{|v|>R}|v|^2\av v^{-\beta}  d v\\
 &\le\tfrac{3}{2}\mathfrak{m}\vartheta+C_\beta N R^{5-\beta}.
\end{align*}
Hence $\mathfrak{m}\le C_\beta N\vartheta^{(3-\beta)/2}$ in this last case.
The three cases combine to give
\begin{equation}\label{eq:app-PP-parameter}
 \mathfrak{m}\vartheta^{-3/2}(1+|u|^2+\vartheta)^{\beta/2}
 \le C_\beta N.
\end{equation}
For every $v$, Gaussian maximization gives
\[
 \av v^\beta e^{-|v-u|^2/(2\vartheta)}
 \le C_\beta(1+|u|^2+\vartheta)^{\beta/2}.
\]
Multiplying this by $\mathfrak{m}(2\pi\vartheta)^{-3/2}$ and using
\eqref{eq:app-PP-parameter} proves \eqref{eq:app-PP-weighted}.

We shall also use the estimate for a sum. The mass, momentum,
and energy of $M_h+M_k$ equal those of $h+k$, so
$M_{M_h+M_k}=M_{h+k}$. Applying
\eqref{eq:app-PP-weighted} to $M_h+M_k$ gives
\begin{equation}\label{eq:app-weighted-sum}
 P_\beta(h+k)
 \le C_\beta N_\beta(M_h+M_k)
 \le C_\beta\bigl(P_\beta(h)+P_\beta(k)\bigr).
\end{equation}
This allows us to estimate the recent BGK source separately from
the older part of the solution.

To control the older part, fix $s\ge2$, $x\in\T^3$, and a delay
$\tau>0$. For each
$y\in[0,2\pi)^3$ and $n\in\mathbb Z^3$, put
$v_n=(x-y+2\pi n)/\tau$. The elementary one-dimensional Gaussian
lattice estimate
\begin{equation}\label{eq:1-d gaussian}
 \sup_{a\in\mathbb R}\sum_{j\in\mathbb Z}
 e^{-\alpha(j-a)^2}\le C(1+\alpha^{-1/2}),\qquad \alpha>0.
\end{equation}
(This can be done by choosing a nearest integer to $a$: there are at most two
integers at each successive distance, so the sum is bounded by
$2+2\sum_{j\ge1}e^{-\alpha(j-1/2)^2}$, which is at most
$C(1+\int_0^\infty e^{-\alpha r^2/4}  d r)$.
The latter equals $C(1+\alpha^{-1/2})$.) 

Applying \eqref{eq:1-d gaussian} in each coordinate gives, uniformly
in the translation of the lattice,
\begin{equation}\label{eq:app-lattice}
 \tau^{-3}\sum_{n\in\mathbb Z^3}
 e^{-|v_n-U|^2/(4T)}
 \le C\tau^{-3}(1+\tau\sqrt T)^3
 \le C\bigl(T^{3/2}+\tau^{-3}\bigr).
\end{equation}
Furthermore,
$\av{v_n}^{2}\le C(1+|U|^2+T)
 (1+|v_n-U|^2/T)$. Absorbing the last factor into half of the
Gaussian exponent in \eqref{eq:app-lattice} yields
\begin{equation}\label{eq:app-weighted-lattice}
 \tau^{-3}\sum_{n\in\mathbb Z^3}
 \av{v_n}^{2}e^{-|v_n-U|^2/(2T)}
 \le C\bigl(T^{3/2}+\tau^{-3}\bigr)
          (1+|U|^2+T).
\end{equation}
Periodic unfolding has Jacobian
$  d v=(2\pi)^3\tau^{-3}  d y$ when $  d y$ is normalized torus
measure. Thus \eqref{eq:app-weighted-lattice}, with
$q_F=\rho(2\pi T)^{-3/2}$, implies
\begin{align}
 &\int_{\R^3}\av v^2 M_F(s,x-\tau v,v)  d v\nonumber\\
 &\quad\le C\int_{\T^3}q_F(s,y)
       \bigl(T(s,y)^{3/2}+\tau^{-3}\bigr)
       (1+|U(s,y)|^2+T(s,y))  d y\nonumber\\
 &\quad\le C(1+\tau^{-3})
       \int_{\T^3}\rho(s,y)
       (1+|U(s,y)|^2+T(s,y))  d y\nonumber\\
 &\quad\le C(1+\tau^{-3}),\quad \tau > 0.
 \label{eq:app-trace}
\end{align}
Here we used $q_FT^{3/2}=(2\pi)^{-3/2}\rho$ and
$q_F\le C\rho$, the latter following from
$T\ge b_0/(3A_\rho)$. For the final line, conservation gives
$\int\rho  d y=1$ and
$\int\rho(|U|^2+3T)  d y=3$. The constant in
\eqref{eq:app-trace} is independent of $s,x,\tau$.
We treat short delays separately because of the factor $\tau^{-3}$.

Fix $\delta\in(0,1)$, to be chosen at the end. Restart the mild
formula at time $2$ and write $F=R_\delta+O_\delta$, where
\begin{align*}
 R_\delta(t,x,v)
 &=\int_0^{\min\{\delta,t-2\}}e^{-\tau}
       M_F(t-\tau,x-\tau v,v)  d\tau,\\
 O_\delta(t,x,v)
 &=e^{-(t-2)}F(2,x-(t-2)v,v)
   +\int_\delta^{t-2}e^{-\tau}
       M_F(t-\tau,x-\tau v,v)  d\tau.
\end{align*}
The last integral is zero when $t-2<\delta$; both terms are
nonnegative. For a finite $\mathfrak{T}>4$, define
$X_{\mathfrak{T}}=\sup_{2\le t\le \mathfrak{T}}N_\beta(F(t))<\infty$, using the
finite-horizon bound of \cite{PP93}. By
\eqref{eq:app-PP-weighted}, the recent part satisfies
\begin{equation}\label{eq:app-recent}
 N_\beta(R_\delta(t,x,\cdot))
 \le C_\beta X_{\mathfrak{T}}\int_0^{\min\{\delta,t-2\}}e^{-\tau}  d\tau
 \le C_\beta\delta X_{\mathfrak{T}},\qquad 2\le t\le \mathfrak{T}.
\end{equation}

The older part has a uniform local second moment. Indeed,
$\beta>5$ makes
$I_\beta=\int_{\R^3}\av v^{2-\beta}  d v$ finite, so the
ballistic term has moment at most $I_\beta N_\beta(F(2))$.
For the source term, use \eqref{eq:app-trace} at the time
$s=t-\tau\ge2$ and then integrate only over $\tau\ge\delta$:
\begin{equation}\label{eq:app-old-moment}
 \sup_{t\ge2,\,x\in\T^3}
 \int_{\R^3}\av v^2O_\delta(t,x,v)  d v
 \le I_\beta N_\beta(F(2))
      +C\int_\delta^\infty e^{-\tau}(1+\tau^{-3})  d\tau
 =:E_\delta<\infty.
\end{equation}
The constant $E_\delta$ may deteriorate as $\delta\downarrow0$,
but is independent of the endpoint ${\color{black}\mathfrak{T}}$. \\

For $t\ge4$, the source integral defining $O_\delta$ contains all
delays in $[1,2]$, since $\delta<1$ and $t-2\ge2$.
Consequently $O_\delta(t,x,v)\ge B_t(x,v)$. Denote its local
macroscopic fields by $\rho_O,U_O,T_O$. By
\eqref{eq:app-bath-moments} and \eqref{eq:app-old-moment},
\[
 b_0\le\rho_O\le E_\delta,
 \qquad |U_O|\le E_\delta/b_0,
 \qquad 3\rho_O T_O
   =\int|v-U_O|^2O_\delta  d v
   \ge b_0(1+|U_O|^2).
\]
Also $3\rho_OT_O\le\int|v|^2O_\delta  d v\le E_\delta$.
In particular,
\begin{equation}\label{eq:app-old-parameters}
 \frac{b_0}{3E_\delta}\le T_O
 \le\frac{E_\delta}{3b_0},
 \qquad \rho_O\in[b_0,E_\delta],\quad
 |U_O|\le E_\delta/b_0.
\end{equation}
These compact parameter bounds give
$P_\beta(O_\delta(t,x,\cdot))\le C_\delta$ for $t\ge4$ by
maximizing its Maxwellian explicitly. 

For $2\le t\le4$ we instead
use $0\le O_\delta\le F$ and the finite-horizon bound:
$P_\beta(O_\delta)\le C_\beta N_\beta(O_\delta)
\le C_\beta\sup_{2\le s\le4}N_\beta(F(s))$.
Enlarging $C_\delta$ gives
\begin{equation}\label{eq:app-old-weighted}
 \sup_{t\ge2,\,x\in\T^3}
 P_\beta(O_\delta(t,x,\cdot))\le C_\delta<\infty.
\end{equation}
The bound \eqref{eq:app-old-moment} concerns the old component only.
No pointwise energy bound for $F$ has been used.

Apply the sum estimate \eqref{eq:app-weighted-sum} at every
$(t,x)$ to $F(t,x,\cdot)=R_\delta(t,x,\cdot)+
O_\delta(t,x,\cdot)$. Equations \eqref{eq:app-recent} and
\eqref{eq:app-old-weighted} imply
\[
 \sup_{\substack{2\le t\le \mathfrak{T} \\x\in \mathbb{T}^3}} P_\beta(F(t,x,\cdot))
 \le C_\beta C_\delta+C_\beta^2\delta X_{\mathfrak{T}}.
\]
The mild formula from time $2$ now gives, since the exponential
kernel integrates to at most $1$,
\[
 X_{\mathfrak{T}}\le N_\beta(F(2))+C_\beta C_\delta
                         +C_\beta^2\delta X_{\mathfrak{T}}.
\]
Choose $\delta>0$ sufficiently small that
$C_\beta^2\delta\le1/2$. Although $C_\delta$ may then be large,
it is fixed and independent of $\mathfrak{T}$. Absorption yields
\begin{equation}\label{eq:app-uniform-weight}
 \sup_{t\ge2}\|\av v^\beta F(t)\|_{L^\infty_{x,v}}
 \le 2\bigl(N_\beta(F(2))+C_\beta C_\delta\bigr)
 =:M_\beta<\infty.
\end{equation}

The uniform weighted bound now gives a pointwise upper bound on the
local energy:
\[
 e_2(t,x)\le M_\beta
       \int_{\R^3}|v|^2\av v^{-\beta}  d v
 =:E_*^{\rm up}<\infty,\qquad t\ge2.
\]
Together with \eqref{eq:app-macro-lower} this gives
$|U(t,x)|\le U^{\rm up}:=\sqrt{E_*^{\rm up}/b_0}$ and
$T(t,x)\le T^{\rm up}:=E_*^{\rm up}/(3b_0)$; its lower bound is
$T_-:=b_0/(3A_\rho)$. Hence, 
\[
 M_F(t,x,v)
 \ge\frac{b_0}{(2\pi T^{\rm up})^{3/2}}
       \exp\!\left(-\frac{(|v|+U^{\rm up})^2}{2T_-}\right)
 \ge a_1e^{-b_1|v|^2},\quad t \geq 2
\]
with positive constants depending only on the established bounds.
For $t\ge3$, the positive source over the last unit of time in the
mild formula therefore gives
$F(t,x,v)\ge(1-e^{-1})a_1e^{-b_1|v|^2}$.
Combining this with \eqref{eq:app-uniform-weight} proves
\eqref{eq:app-envelopes}; we may take $t_1=4$.

Theorem~\ref{thm:bgk}, applied after $t_1$, gives
$\HH(t)\le C e^{-\lambda(t-t_1)}$, so $\HH(t)\to0$.
The uniform weighted bound also gives
\[
 \|\av v^\beta(F(t)-\mu)\|_\infty
 \le M_\beta+\|\av v^\beta\mu\|_\infty=:M_*
 \quad(t\ge t_1).
\]
The Gaussian lower bound in \eqref{eq:app-envelopes} supplies
a restart free-streaming density constant
\[
 C_{\rm fs}:=a\int_{\R^3}e^{-b|v|^2}  d v>0:
 \qquad
 \inf_{\tau\ge0,\,x\in\T^3}
 \int_{\R^3}F(s,x-\tau v,v)  d v\ge C_{\rm fs}
 \quad(s\ge t_1).
\]
Choose a finite time $s\ge t_1$ at which
$\HH(s)$ is below the small-entropy threshold
$\varepsilon(M_*,C_{\rm fs})$ of \cite[Theorem 1.1]{BKLY26}.
Such a time exists by the preceding entropy decay. The collision
frequency $1$ corresponds to both exponents being zero in the family
$\nu=\rho^aT^b$ of \cite{BKLY26}, and $\beta>10$ meets its
weighted exponent requirement. Conservation gives its equilibrium
parameters; uniqueness of the time-translated mild solution
identifies the restart with the solution covered by that theorem.
It follows that, for some $C_\infty,\kappa>0$,
\[
 \|\av v^\beta(F(t)-\mu)\|_\infty
 \le C_\infty e^{-\kappa(t-s)},\qquad t\ge s.
\]

We finish by passing from weighted $L^\infty$ decay to entropy decay.
Write
$\Phi(F\mid\mu)=F\log(F/\mu)-F+\mu$. On $F<\mu$,
$F\log(F/\mu)\le0$, hence $\Phi(F\mid\mu)\le\mu-F$.
On $F\ge\mu$, the inequality $r-1-\log r\ge0$ for
$r=F/\mu\ge1$ gives
\[
 \Phi(F\mid\mu)\le(F-\mu)\log(F/\mu).
\]
Since $F\le M_\beta\av v^{-\beta}$ and
$\mu=(2\pi)^{-3/2}e^{-|v|^2/2}$, on that same set
$0\le\log(F/\mu)\le C(1+|v|^2)$.
Both cases therefore satisfy
\[
 0\le\Phi(F\mid\mu)
 \le C(1+|v|^2)|F-\mu|.
\]
Integrating and using $\beta>5$ yields
\[
 \HH(t)\le C\left(\int_{\R^3}
               (1+|v|^2)\av v^{-\beta}  d v\right)
              \|\av v^\beta(F(t)-\mu)\|_\infty
 \le C_H e^{-\kappa(t-s)}.
\]
Taking $t_{\rm ent}=s$ proves \eqref{eq:app-weighted-decay}.
\end{proof}

\begin{remark}[Scope of the application]
\label{rem:app-scope}
Theorem~\ref{thm:bgk} in the present draft already gives
\emph{exponential}, rather than merely polynomial, entropy decay
once \eqref{eq:app-envelopes} is available. Only $\HH(t)\to0$ is
needed to enter the small-entropy regime of \cite{BKLY26}.
The additional assumption \eqref{eq:app-rho-upper} is essential to
the argument above: the time-dependent bounds in \cite{PP93}
alone cannot be substituted for a uniform density bound. \\
\end{remark}

\bigskip

\noindent{\bf Data availability:} No data was used for the research described in the article.
\newline

\noindent{\bf Conflict of interest:} The authors declare that they have no conflict of interest.\newline

\noindent{\bf Disclaimer on the use of AI tools} During the preparation of this manuscript, the author utilized generative AI (ChatGPT) exclusively for English language editing, identifying relevant references, cross-checking the validity of mathematical derivations, and refining the mathematical rigor of Lemma \ref{lem:BE-chain} and Lemma \ref{aux:villani}. All core mathematical ideas, formulations, and results originated entirely from the authors. The authors take full responsibility for the contents and ultimate integrity of this work. \newline

\noindent{\bf Acknowledgement} D. Lee is supported by the National Research Foundation of Korea(NRF) grant funded by the Korea government(MSIT)(No.RS-2023-00212304 and No.RS-2023-00219980). \newline

\bigskip

\appendix
\section{An auxiliary-kernel comparison for polynomial tails}
\label{sec:aux}

Section 6.1 of \cite{HJJLP26} uses Villani's entropy-production
inequality \cite{V03} and H\"older's inequality. With H\"older
exponent $3/2$, that argument gives microscopic power $3/2$ for a
polynomial tail exponent greater than eleven. The decay rate is not optimal, of course. We use a different
comparison: after fixing the constant in Villani's inequality on
the envelope class, we separate small, intermediate, and large
relative velocities. 
Under the uniform angular lower bound $b\ge c_b>0$ almost
everywhere, the resulting power depends on the tail exponent and
the comparison requires only $m>7$ for the velocity weight. 

For a positive velocity distribution $f$, write and recall
\[
 u:=|v-v_*|,\qquad
 \Psi(A,B)=(A-B)\log(A/B),\qquad
 h=H_v(f\mid M_f),
\]
and define
\begin{align}
 D_v(f)
 &:=\frac14\int_{\R^3\times\R^3\times\mathbb S^2}
 u^\gamma b(\widehat{v-v_*}\cdot\sigma)
 \Psi(f'f_*',ff_*)\,  d\sigma\,  d v_*\,  d v,
 \label{aux:physical}\\
 D_{\mathrm{aux}}(f)
 &:=\frac14\int_{\R^3\times\R^3\times\mathbb S^2}
 (1+u^2)\Psi(f'f_*',ff_*)\,  d\sigma\,  d v_*\,  d v.
 \label{aux:production}
\end{align}
The postcollisional velocities and angular convention are the same as
in the Boltzmann equation.

\begin{lemma}
\label{aux:villani}
Fix $A,a,M>0$ and $m>7$. For every positive $f$ such that
\begin{equation}\label{aux:envelopes}
 A e^{-a|v|^2}\le f(v)\le M\av v^{-m},
\end{equation}
there is a constant $\lambda_{\mathrm{aux}}>0$, depending only on
$A,a,M,m$ and the fixed collision convention, for which
\begin{equation}\label{aux:linear}
 D_{\mathrm{aux}}(f)\ge
 \lambda_{\mathrm{aux}}H_v(f\mid M_f).
\end{equation}
\end{lemma}
\begin{proof}
Write $\rho,U,T$ for the mass, mean velocity, and temperature of
$f$. The two envelopes imply
\[
 \rho_-:=A\int e^{-a|v|^2}  d v\le\rho\le
 \rho_+:=M\int\av v^{-m}  d v,\qquad
 E_+:=M\int|v|^2\av v^{-m}  d v<\infty,
\]
and therefore
\[
 T\le T_+:=\frac{E_+}{3\rho_-},\qquad
 T\ge T_-:=\frac{A}{3\rho_+}
                  \int|v|^2e^{-a|v|^2}  d v>0.
\]
In the last bound, the cross term in
$\int |v-U|^2e^{-a|v|^2}  d v$ vanishes by symmetry.

As in (9) of \cite{V03}, normalize by 
\[
	\widehat f(w)=T^{3/2}\rho^{-1}f(U+\sqrt T\,w).
\]
It has mass one, zero mean, and energy three. {\color{black}\hypersetup{linkcolor=black,citecolor=black}
Define its
minimum directional temperature by $\Theta(\widehat f)$ as below.
The Gaussian lower bound gives
}
\begin{align*}
 \Theta(\widehat f)
 &:=\inf_{|e|=1}\int(w\cdot e)^2\widehat f(w)  d w\\
 &=\inf_{|e|=1}\frac{1}{\rho T}
       \int((v-U)\cdot e)^2f(v)  d v\\
 &\ge\frac{A}{\rho_+T_+}
       \inf_{|e|=1}\int((v-U)\cdot e)^2e^{-a|v|^2}  d v\\
 &\ge\frac{A}{2a\rho_+T_+}
                    \left(\frac{\pi}{a}\right)^{3/2}
 =:\Theta_->0,
\end{align*}
since the last integral equals
$(\pi/a)^{3/2}((2a)^{-1}+(U\cdot e)^2)$.
{\color{black}\hypersetup{linkcolor=black,citecolor=black}
Let $\lambda_1\le\lambda_2\le\lambda_3$ be the eigenvalues of
$\int w\otimes w\,\widehat f(w)\,d w$. Then
$\lambda_1=\Theta(\widehat f)$ and
$\lambda_1+\lambda_2+\lambda_3=3$. Villani's inequality
\cite[Theorem 2.1]{V03}, for the kernel $1+|w-w_*|^2$ with
constant angular factor one, gives
\[
 D_{\mathrm{aux}}(\widehat f)
 \ge c_{d}(3-\lambda_3)H_v(\widehat f\mid\mu),
\]
where $c_{d}>0$ is dimensional and $\lambda_3$ is the maximum
directional temperature. Since
\[
 3-\lambda_3=\lambda_1+\lambda_2
 \ge2\lambda_1=2\Theta(\widehat f),
\]
we obtain
\[
 D_{\mathrm{aux}}(\widehat f)
 \ge 2c_d\Theta(\widehat f)H_v(\widehat f\mid\mu).
\]
}

Changes of variables and the collision invariances yield
\[
 H_v(f\mid M_f)=\rho H_v(\widehat f\mid\mu),\qquad
 D_{\mathrm{aux}}(f)
 =\rho^2\bigl(D_{0,1}(\widehat f)+T D_{2,1}(\widehat f)\bigr)
 \ge\rho^2\min\{1,T\}D_{\mathrm{aux}}(\widehat f),
\]
where 
\[
D_{j,1}(g)
:=\frac14
\int_{\mathbb R^3}\int_{\mathbb R^3}\int_{\mathbb S^2}
|w-w_*|^j
\Psi\!\left(g(w')g(w_*'),\,g(w)g(w_*)\right)
\,d\sigma\,dw_*\,dw,\ \ j=0,2.
\]
Consequently \eqref{aux:linear} holds with
\[
 \lambda_{\mathrm{aux}}
 =c_{d}\rho_-\min\{1,T_-\}\Theta_->0.
\]
\end{proof}

\begin{proposition} 
\label{aux:proposition}
Let $0\le\gamma\le1$, let $b\ge c_b>0$ almost everywhere, and assume
that $b$ is integrable on the sphere. Let $f$ satisfy
\eqref{aux:envelopes} with $m>7$. Then there exists
$c_{\mathrm{mic}}>0$, depending only on $A,a,M,m,\gamma,{\color{black}c_b}$
and the collision convention, such that
\begin{equation}
 D_v(f)\ge c_{\mathrm{mic}}
 H_v(f\mid M_f)^{1+\varepsilon_*},
 \qquad
 \varepsilon_*=\max\left\{\frac\gamma3,
                           \frac{2-\gamma}{m-7}\right\}.
 \label{aux:power}
\end{equation}
Consequently, a Boltzmann solution in this uniform class satisfying
the hypotheses of Theorem~\ref{thm:boltzmann} obeys
\begin{equation}
 H(F(t)\mid\mu)\le C(1+t-t_*)^{-r},\qquad
 r=\left(\max\left\{\frac\gamma3,
                         \frac{2-\gamma}{m-7}\right\}\right)^{-1}.
 \label{aux:rate}
\end{equation}
For hard spheres, $\gamma=1$, this gives $r=\min\{3,m-7\}$.
\end{proposition}

\begin{proof}
The envelopes imply $|\log f(v)|\le C\av{v}^2$. Put
$S=1+|v|^2+|v_*|^2$. Energy conservation gives
\begin{equation}
 \left|\log\frac{f'f_*'}{ff_*}\right|\le CS,
 \qquad
 \Psi(f'f_*',ff_*)\le C(f'f_*'+ff_*)S.
 \label{aux:integrand}
\end{equation}
First, \eqref{aux:production} is finite. Since
$(1+u^2)S\le2S^2$, the collision change of variables bounds it by a
constant times
$\iint ff_*S^2\,  d v_*  d v$, and $m>7$ guarantees a finite fourth
moment. The same bounds and the uniform macroscopic bounds implied by
\eqref{aux:envelopes} give
\begin{equation}
 0\le h\le h_{\max}<\infty
 \label{aux:hmax}
\end{equation}
uniformly over the class.

Fix $0<\delta\le1$ and $R\ge2$. Split
$D_{\mathrm{aux}}$ into the regions $u<\delta$,
$\delta\le u\le R$, and $u>R$. The standard pre/postcollisional
symmetry preserves $u$, $S$, and the full collision measure in
\eqref{aux:production}. Thus the gain and loss terms in the upper bound
\eqref{aux:integrand} have equal integrals over each of these regions. 

\textit{Region 1} On $u<\delta$, one has $S\le C\av{v}^2$ and $f_*\le M$. Hence this
part, denoted by $E_\delta$, satisfies
\begin{align}
 E_\delta
 &\le C\iint_{|v-v_*|<\delta}ff_*S\,  d v_*  d v
 \le C\delta^3\int\av{v}^2f(v)\,  d v
 \le C_{\rm small}\delta^3.
 \label{aux:small}
\end{align}
Here $C_{\rm small}$ depends only on the envelopes, which control the second
moment.

\textit{Region 2} For the part $E_R$ with $u>R$, use
\[
 \{u>R\}\subset
 \{|v|>R/2\}\cup\{|v_*|>R/2\},
 \qquad S^2\le2(\av{v}^4+\av{v_*}^4).
\]
Writing $\mathcal M_j(f)=\int\av{v}^j f(v)\,  d v$, we obtain
\begin{align*}
 E_R\le C\bigg[
 \mathcal M_0(f)\int_{|v|>R/2}\av{v}^4f(v)\,  d v
 +\mathcal M_4(f)\int_{|v|>R/2}f(v)\,  d v
 \bigg].
\end{align*}
For $j=0,4$, the upper envelope and spherical integration imply
\[
 \int_{|v|>R/2}\av{v}^j f(v)\,  d v
 \le C_{M,m}\frac{1}{R^{m-7}}.
\]
It follows that
\begin{equation}
 E_R\le C_{\rm large} R^{-(m-7)}.
 \label{aux:large}
\end{equation}
Here $C_{\rm large}$ depends only on $A,a,M,m$ and the fixed collision convention, and is independent of $R$.

\textit{Region 3} On the intermediate region,
\[
 \frac{1+u^2}{u^\gamma}
 =u^{-\gamma}+u^{2-\gamma}
 \le\delta^{-\gamma}+R^{2-\gamma}.
\]
Using $b\ge {\color{black}c_b}$ there, and then \eqref{aux:linear}, yields
\begin{equation}
 \lambda_{\mathrm{aux}}h
 \le c_b^{-1}(\delta^{-\gamma}+R^{2-\gamma})D_v(f)
      +C_{\rm small}\delta^3+C_{\rm large} R^{-(m-7)}.
 \label{aux:comparison}
\end{equation}
Enlarge $C_{\rm small},C_{\rm large}$ to positive constants if necessary. Choose a fixed
$h_0\in(0,1]$ sufficiently small that, for $0<h\le h_0$,
\begin{equation}
 \delta=\left(\frac{\lambda_{\mathrm{aux}}h}{4C_{\rm small}}\right)^{1/3}
 \le1,\qquad
 R=\left(\frac{4C_{\rm large}}{\lambda_{\mathrm{aux}}h}\right)^{1/(m-7)}
 \ge2.
 \label{aux:cutoffs}
\end{equation}
The last two terms in \eqref{aux:comparison} then sum to
$\lambda_{\mathrm{aux}}h/2$. Moreover,
\[
 \delta^{-\gamma}+R^{2-\gamma}
 \le C\left(h^{-\gamma/3}
            +h^{-(2-\gamma)/(m-7)}\right)
 \le C' h^{-\varepsilon_*}.
\]
We have proved \eqref{aux:power} when $0<h\le h_0$.
For $h\ge h_0$, keep the cutoffs chosen at $h_0$.
The error in \eqref{aux:comparison} is at most
$\lambda_{\mathrm{aux}}h_0/2\le\lambda_{\mathrm{aux}}h/2$, so
$D_v(f)\ge c_0h$. In view of \eqref{aux:hmax}, this implies
$D_v(f)\ge c_0\max\{1,h_{\max}\}^{-\varepsilon_*}h^{1+\varepsilon_*}$.
The case $h=0$ is immediate. 
\end{proof}



\end{document}